\documentclass{amsart}
\usepackage{amssymb, amsmath, amsfonts, amscd}

\input xypic

\theoremstyle{plain}
\newtheorem{theorem}{Theorem}[section]
\newtheorem{corollary}[theorem]{Corollary}
\newtheorem{lemma}[theorem]{Lemma}
\newtheorem{proposition}[theorem]{Proposition}

\newtheorem{example}[theorem]{Example}

\theoremstyle{definition}
\newtheorem{definition}[theorem]{Definition}

\theoremstyle{remark}

\numberwithin{equation}{theorem}

\newcommand{\Hom}{\operatorname{Hom} }

\newcommand{\End}{\operatorname{End} }

\renewcommand{\H}{\operatorname{H}}
\newcommand{\Sym}{\operatorname{Sym} }
\newcommand{\SLO}{\operatorname{SL} }
\newcommand{\SL }{\operatorname{SL}_{2,A} }

\newcommand{\Spec}{\operatorname{Spec} }

\newcommand{\C}{\mathbf{C} }
\newcommand{\Q}{\mathbf{Q} }

\renewcommand{\sl}{\mathfrak{sl}_{2,A} } 

\newcommand{\glo}{\mathfrak{gl} }
\newcommand{\gl}{\mathfrak{gl}_{2,A} }

\newcommand{\Ad}{\operatorname{Ad} }
\newcommand{\GL}{\operatorname{GL}_{2,A} }
\newcommand{\GLO}{\operatorname{GL}}

\newcommand{\Lie}{\operatorname{Lie} }
\newcommand{\Dist}{\operatorname{Dist} }
\newcommand{\U}{\operatorname{U}}

\newcommand{\z}{x_{21}}
\newcommand{\w}{x_{22} }
\newcommand{\y}{x_{12} } 
\newcommand{\x}{x_{11} }

\newcommand{\Z}{z_{21}}
\newcommand{\W}{z_{22} }
\newcommand{\Y}{z_{12} } 
\newcommand{\X}{z_{11} }

\newcommand{\alg}{\underline{A-alg}}
\newcommand{\grp}{\underline{Groups}}
\newcommand{\ZZ}{\mathbf{Z}} 

\renewcommand{\C}{\mathbf{C}}
\newcommand{\gm}{\mathbb{G}_m} 
\newcommand{\ga}{\mathbb{G}_a} 
\newcommand{\sym}{\operatorname{S}}

\begin{document}

\title{A refined Weyl character formula for comodules on $\GL$}

\author{Helge \"{O}ystein  Maakestad }

\address{Sætre i Hurum, Buskerud, Norway} 


\date{April 2024} 

\begin{abstract} Let $A$ be any commutative unital ring and let $\GL$ be the general linear group on $A$ of rank $2$. I study the representation theory of $\GL$ and the symmetric powers $\Sym^d(V)$, where $(V, \Delta)$ is the standard right comodule 
on $\GL$. I prove a \emph{refined Weyl character formula} for $\Sym^d(V)$. There is for any integer $d \geq 1$  a (canonical) refined weight space decomposition $\Sym^d(V)  \cong  \oplus_{i } \Sym^d(V)^i$ where each direct summand  $\Sym^d(V)^i$ is a comodule on $N \subseteq \GL$. Here $N$ is the schematic normalizer of the diagonal torus $T \subseteq \GL$. I prove a character formula for the direct summands of $\Sym^d(V)$ for any integer $d \geq 1$. This refined Weyl character formula implies the classical Weyl character formula. As a Corollary I prove a refined Weyl character formula for the pull back $\Sym^d(V \otimes K)$ as a comodule on $\operatorname{GL}_{2,K}$ where $K$ is any field.   I also prove that the direct summand $\Sym^d(V)^i$ is an irreducible $N$-comodule over any field. I moreover calculate explicit examples involving the symmetric powers, symmetric tensors and their duals. 
\end{abstract}

\maketitle
\tableofcontents

\section{Introduction} 
In a recent paper (see \cite{maa1}) I studied the notion "good filtration" for torsion free comodules on $\SLO_{2,\ZZ}$ and I proved the existence of a "virtual Clebsch-Gordan formula" for the tensor product $\Sym^m(V) \otimes \Sym^n(V)$ of the symmetric powers of the standard right comodule $(V, \Delta)$ on the group scheme $\SL$ where $A$ is an arbitrary commutative unital ring and $d\geq 1$ is an integer. This gives an infinite set of non trivial examples of torsion free finite rank comodules on $\SLO_{2,\ZZ}$ equipped with a good filtration. In general one wants to give a definition of the notion "good filtration" and to classify torsion free finite rank comodules with a good filtration.   In \cite{maa1} I found examples of non isomorphic comodules on $\SLO_{2,\ZZ}$ that become isomorphic when we "restrict to the generic fiber". There is a canonical morphism of group schemes $\pi: \SLO_{2,\mathbf{Q}} \rightarrow \SLO_{2,\ZZ}$ and I proved existence of non isomorphic finite rank  torsion free comodules $W_1,W_2$ with $\pi^*W_1 \cong \pi^*W_2$ an isomorphism of comodules on $\SLO_{2,\mathbf{Q}}$.  In order to define the notion "good filtration" for comodules on $\SLO_{2,\ZZ}$, we must for each integer $d\geq 1$ classify comodules $W$ on $\SLO_{2,\ZZ}$ with the property  that the pull back to $\SLO_{2,\mathbf{Q}}$ becomes isomorphic to $\Sym^d(V)$ where  $V:=\mathbf{Q}^2$ is the "standard right comodule" of $\SLO_{2,\mathbf{Q}}$. There is by \cite{maa1} no highest weight theory or "complete reducibility property" for torsion  free comodules on $\SLO_{2,\ZZ}$ and $\GLO_{2, \ZZ}$, and for this reason we need new invariants and new methods in the study of such comodules and the aim of this paper is to introduce a new structure on a class of torsion free comodules on $\SLO_{2,\ZZ}$. 

In this paper I prove the following: Let $A$ be any commutative unital ring and let $G:=\GL$ or $\SL$. I define the \emph{schematic normalizer} $N  \subseteq G$, which is a closed subgroup scheme of $G$ containing the diagonal torus $T$ as a closed normal subgroup scheme. The schematic normalizer is isomorphic to the functor normalizer $N_G(T)$ as defined in \cite {jantzen}.   Let $\sigma, q: N \times T \rightarrow N$ be the action and projection map. 


\begin{theorem} \label{quoGL} Let $T$ be the diagonal torus with $N \subseteq G$ the schematic normalizer. There  is a faithfully flat map of affine group schemes $\pi: N \rightarrow N/T$ with $T \cong ker(\pi)$ giving an exact sequence of group schemes

\[ e \rightarrow T \rightarrow N \rightarrow^{\pi} N/T \rightarrow e .\]

For any map $f: N \rightarrow Z$ to an affine scheme $Z$ with $ f \circ \sigma = f \circ q$ there is a unique map $f^*: N/T \rightarrow Z$ with $f^* \circ \pi = f$. Hence the affine group scheme $N/T$ is uniquely determined up to isomorphism.
If $G:=\GL$ it follows the sequence splits and there is an isomorphism of group schemes $N \cong T \times^{\sigma} N/T$  where the product on the right is the semi direct product.

\end{theorem}

The quotient in \ref{quoGL} is proved to exist in \cite{jantzen}, \cite{DG} and \cite{waterhouse} over an arbitrary Dedekind domain, but the ring $A$ in Theorem \ref{quoGL} may be an arbitrary commutative unital ring. I give an explicit construction
of the representing $A$-Hopf algebra $k[N/T]$ of $N/T$ and the morphism $\pi$. If $W_G$ is the Weyl group of $\GL$ it follows there is an isomorphism of group schemes

\[ N \cong T \times^{\sigma} \Spec(k[W_G]) \]

where $\Spec(k[W_G])$ is the constant group scheme of the Weyl group $W_G$.

Then I prove the existence of a \emph{refined weight space decomposition} of any symmetric power $W:=\Sym^d(V)$ where $d \geq 1$ is an integer and where $(V, \Delta)$ is the standard right comodule on $\GL$. There is a canonical decomposition $W \cong \oplus_{\lambda \in \Lambda} W_{\lambda}$ where each summand $W_{\lambda}$ is a comodule on the schematic normalizer $N \subseteq \GL$ of the diagonal torus $T \subseteq \GL$. Then I prove a \emph{refined Weyl character formula} for each summand $W_{\lambda}$ of $W$. The refined character formula implies the classical Weyl character formula (this is Theorem \ref{refined} in the paper):

\begin{theorem} \label{mainintro}  Let $d \geq 1$  be an integer and let $\Sym^d(V) \cong \oplus_i \Sym^d(V)^i$ be the refined weight space decomposition of $\Sym^d(V)$.  There is  for each integer $i \geq 0$ a $k[N]$-comodule structure

\[ \Delta_{N,d,i}: \Sym^d(V)^{i} \rightarrow k[N]\otimes \Sym^d(V)^{i}, \]

and a direct sum decomposition $\Sym^d(V) \cong \oplus_{ i } \Sym^d(V)^{i}$ of $k[N]$-comodules.
The character of $\Sym^d(V)^i$ is the following: If $d:=2k+1 \geq 1$ is an odd integer it follows 

\[ Char(\Sym^d(V)^i)=x_2^{d-i}x_1^i + x_2^i x_1^{d-i}\]

 for $i=0,..,k$. If $d:=2k \geq 2$  the following holds:

\[ Char(\Sym^d(V)^i)= x_1^{d-i}x_2^i +x_1^i x_2^{d-i} \]

if $i=0,..,k-1$, and  

\[  Char(\Sym^d(V)^{k})=x_1^kx_2^k.\]

\end{theorem}

In Example \ref{finitefield} we relate the refined weight space decomposition to irreducible finite dimensional comodules over a field $k$ of positive charcacteristic.  If $G:=\SLO_{n,k}$ and if $S_{\lambda}(V)$ is a Weyl module,
there is a unique finite dimensional  irreducible comodule $V(\lambda) \subsetneq S_{\lambda}(V)$. The irreducible comodule $V(\lambda)$ has a relation with the refined weight space decomposition for $S_{\lambda}(V)$.

In Theorem \ref{irreducibleN} I prove that the direct summand $\Sym^d(V)^i$ is an irreducible comodule on $N$ over any field $k$. In general let  $V$ be a finite dimensional vector space over a field $k$ of characteristic $p>0$ and let $S_{\lambda}(V)$ be the Weyl module of $V$ corresponding to the partition $\lambda$. If $G:=\GLO(V)$ there is in general a unique irreducible sub-$G$-comodule 

\[ V(\lambda) \subseteq S_{\lambda}(V) .\]

There are refined weight space decompositions $S_{\lambda}(V) \cong \oplus_{i \in I} S_{\lambda}(V)^i$ and $V(\lambda) \cong \oplus_{j \in J} V(\lambda)^j$ and there are relations between these two decompositions. In the case when $V \cong k^2$ we may express the irreducible comodule $V(\lambda)$ in terms of the refined weight space decomposition for $\Sym^d(V)$. Hence the construction may have applications in the study of irreducible $G$-comodules over any field $k$.

As an application of the methods introduced in the paper I prove the following Theorems (see Theorem \ref{example}). Let $G:=\GL$ where $A$ is any commutative ring and  let $(V, \Delta)$ be the standard left comodule on $G$.
Let  $\Sym^2(V), \sym_2(V), \Sym^2(V)^*$ and  $\sym_2(V)^*$ be the symmetric power, the symmetric tensor and the dual of the symmetric power and symmetric tensor of $V$. The symmetric tensor is the invariants $\sym_2(V) \subseteq V \otimes_A V$ under the canonical action of the symmetric group on two elements.

\begin{theorem}  \label{exampleintro} Let $G:=\GL$ with $A$ any commutative ring, and let $N \subseteq G$ be the schematic normalizer of the diagonal torus. There is an isomorphism of $k[N]$-comodules (and $k[T]$-comodules) $\Sym^2(V) \cong \sym_2(V)$.
Let $G:=\GLO_{2,\ZZ}$ with $\ZZ$ the ring of integers. It follows the $k[G]$-comodules $\Sym^2(V)$ and $\sym_2(V)$ are not isomorphic.
\end{theorem} 

When taking the dual of the isomorphism $\Sym^2(V) \cong \sym_2(V)$ from Theorem \ref{exampleintro} we get an isomorphism of $k[N]$-comodules $\Sym^2(V)^* \cong \sym_2(V)^*$.

It is hoped the results in this paper  may be interesting in the study of torsion free comodules on $\SL$ and $\GL$ in the future. One problem one wants to study is the problem of classifying $\SLO_{2,\ZZ}$    comodules $W$ whose base change to 
$\SLO_{2, \Q}$ are isomorphic to $\Sym^d(\Q^2)$ where $\Q^2$ is the standard comodule on $\SLO_{2,\Q}$. Such a classification may have interest in the problem of giving a definition of the notion "good filtration" for torsion free finite rank comodules on 
$\SLO_{2,\ZZ}$.
Such a notion may have applications in the study of the representation theory of $\SLO_{n,k}$ and $\GLO_{n,k}$ over an arbitrary field $k$.


\section{The adjoint representation of $\GL$ and $\SL$}

Let $A$ be any commutative unital ring and let $x_{11}, x_{12}, x_{21}, x_{22}$ be independent varaibles over $A$. Let $D:=x_{11}x_{22}-x_{12}x_{21}$ be the "determinant" and let $F(x_{ij}):=D-1 \in A[x_{ij}]$. Let $k[S]:=A[x_{ij}]/(F)$ and let $k[G]:=A[x_{ij}, 1/D]$. Let $\SL:=\Spec(k[S])$ and let $\GL:=\Spec(k[G])$.
Let $\alg$ be the category of commutative $A$-algebras and let $\grp$ be the category of groups. Define for $R=k[S], k[G]$ the functor

\[ h_R: \alg \rightarrow \grp \]

by $h_R(B):=\Hom_{A-alg}(R,B)$. It follows $h_R(B)$ is the group of $2 \times 2$ matrices $g:=(a_{ij})$ with coefficients $a_{ij}$ in $B$. If $R:=k[S]$ It follows  $det(g)=1$ and if $R=k[G]$ it follows $det(g)$ is a unit in $B$.
The affine schemes $\GL, \SL$ are affine group scheme in the sense of the book of Demazure/Gabriel and Jantzen and the rings $k[G], k[S]$ are A-Hopf algebras. Define the following structure as a Hopf algebra on $k[G]$ and $k[S]$:

\[ \Delta: R \rightarrow R \otimes_A R \]

by

\[ \Delta(x_{11}):= x_{11} \otimes x_{11} + x_{12} \otimes x_{21} ,\]

\[ \Delta(x_{12}):= x_{11} \otimes x_{12} + x_{12} \otimes x_{22} ,\]

\[ \Delta(x_{21}):= x_{21} \otimes x_{11} + x_{22} \otimes x_{21}, \]

and

\[ \Delta(x_{22}):= x_{21} \otimes x_{12} + x_{22} \otimes x_{22} .\]

Define the coinversion map $S_1: k[G] \rightarrow k[G]$ by

\[ S_1(x_{11}):=x_{22}/D,  S_1(x_{12}):=-x_{12}/D,  S_1(x_{21}):=-x_{21}/D,  S_1(x_{22}):=x_{11}/D \]

and define the coinversion map $S_2: k[S] \rightarrow k[S]$ as follows:

\[ S_2(x_{11}):=x_{22},  S_2(x_{12}):=-x_{12},  S_2(x_{21}):=-x_{21},  S_2(x_{22}):=x_{11}. \]

Define the counit $\epsilon: R \rightarrow A$  for $R=k[G], k[S]$ as follows: 

\[ \epsilon(x_{11}):= 1 , \epsilon(x_{12}):= 0,    \epsilon(x_{21}):=0,    \epsilon(x_{22}):= 1   .\]

It follows the 4-tuples $(k[G], \Delta, S_1, \epsilon)$ and  $(k[S], \Delta, S_2, \epsilon)$ Are A-Hopf algebras in the sense of Jantzens book \cite{jantzen}. If $B$ is any commutative $A$-algebra it follows $h_{k[S]}(B)$ is the group of $2 \times 2$-matrices with coefficients in $B$ with determinant equal to one, and the above defined $A$-Hopf algebra structure on $k[S]$ induce matrix multiplication  as group structure on each set $h_{k[S]}(B)$.  Similarly it follows the set $h_{k[G]}(B)$ is the set
of $2\times 2$ matrices with coefficients in $B$ and determinant a unit in $B$. The comultiplication $\Delta$ induce matrix multiplication as group operation on $h_{k[G]}(B)$ for any $B$.

In this section we define the adjoint representation of $\GL, \SL$  on their Lie algebras using the functor of points $h_G(-)$. 


Define the following functor

\[ \Lie(G)^F(-) : \alg \rightarrow \grp \]

where

\[ \Lie(G)^F(B):=\{ g \in h_G(B[\epsilon]): p_*(g)=Id \in h_G(B)\}   .\]

It follows $\Lie(G)^F(B)$ is the set of all $2\times 2$-matrices with coefficients in $B$. Hence $\Lie(G)^F(B) \cong \glo_{2,B} \cong \gl  \otimes_A B$.

Define the following functor 

\[ (\gl)_a(-): \alg \rightarrow \grp \]

by

\[ (\gl)_a(B):= (\gl \otimes_A B, +) ,\]

where $(\gl \otimes_A B,+)$ is the underlying abelian group  of the $B$-module $\gl\otimes_A B$. For any map of $A$-algebras $\phi: B \rightarrow B_1$ we get an induced map

\[ \phi_a: (\gl)_a(B) \rightarrow (\gl)_a(B_1) \]

defined by $\phi_a(x \otimes b):= x \otimes \phi(b)$, hence the above definition gives a functor as claimed.    There is an isomorphism of functors $\Lie(G)^F(-) \cong (\gl)_a(-)$.    

We want to define an action of $h_G(-)$ on $(\gl)_a(-)$.
Define for any $B \in \alg$ the following action:

\[ \sigma_B(-,-): h_G(B) \times (\gl)_a(B) \rightarrow (\gl)_a(B) \]

by

\[ \sigma_B(g, z):= g z g^{-1}. \]

We have defined a natural transformation of functors

\[ \Ad^F:  h_G(-) \rightarrow h_{\GLO(\gl)}(-) \]

called the \emph{(left) adjoint representation} of $G$ on $\gl$.  There is similarly a right adjoint representation. We want to calculate the corresponding comodule map of the adjoint representation $\Ad^F$.  We get a map

\[ \sigma_{k[G]}(-,-): h_G(k[G]) \times (\gl)_a(k[G]) \rightarrow (\gl)_a(k[G]) \]

and a map

\[  \rho: (\gl)_a(k[G]) \rightarrow (\gl)_a(k[G]) \]

defined by 

\[ \rho( z \otimes f):= \sigma(id_{k[G]}, z \otimes f) \in (\gl)_a(k[G]):=\gl \otimes_A k[G].\]

Define the map 

\[ \Ad: \gl  \rightarrow \gl  \otimes_A k[G] \]

by

\[ \Ad(z):= \rho( z \otimes 1_{k[G]}) \in \gl \otimes_A k[G]. \]

The map $id_{k[G]} \in h_G(k[G])$ corresponds to the following matrix:

\[ 
Id_{k[G]}=
\begin{pmatrix}   \x  &  \y \\
                             \z  & \w 
\end{pmatrix}.
\]

The matrix $id_{k[G]}^{-1}$ corresponds to the following matrix:

\[ 
Id_{k[G]}^{-1}=
\begin{pmatrix}   \w/D  &  -\y/D \\
                             -\z/D  & \x/D 
\end{pmatrix}.
\]

Define the following matrices in $\gl $:

\[ 
e_{11} =
\begin{pmatrix}   1  &  0 \\
                             0   & 0  
\end{pmatrix},
\]

\[ 
e_{12} =
\begin{pmatrix}   0  &  1 \\
                             0   & 0  
\end{pmatrix},
\]

\[ 
e_{21} =
\begin{pmatrix}   0  &  0 \\
                             1   & 0  
\end{pmatrix},
\]

and

\[ 
e_{22} =
\begin{pmatrix}   0  &  0 \\
                             0   & 1  
\end{pmatrix}.
\]

It follows $\gl  \cong A\{e_{ij} \}$ is a free $A$-module of rank $4$ on the basis $B:=\{ e_{ij} \}$. We want to calculate the element $\Ad(e_{ij})$ for all $i,j$.

An explicit  calculation gives the following formulas:

\[ \Ad(e_{11})= e_{11} \otimes (\x\w/D)  -  e_{12} \otimes (\x\y/D) +   e_{21} \otimes (\z\w/D)  -  e_{22} \otimes (\y\z/D) , \]

\[ \Ad(e_{12})= - e_{11} \otimes (\x\z/D)  +  e_{12} \otimes (\x^2/D) -   e_{21} \otimes (\z^2/D)  +  e_{22} \otimes (\x\z/D),  \]

\[ \Ad(e_{21})= e_{11} \otimes (\y\w/D)  -  e_{12} \otimes (\y^2/D) +   e_{21} \otimes (\w^2/D)  -  e_{22} \otimes (\y\w/D)  ,\]

and

\[ \Ad(e_{22})= -  e_{11} \otimes (\y\z/D)  +   e_{12} \otimes (\x\y/D)  -    e_{21} \otimes (\z\w/D)  +   e_{22} \otimes (\x\w/D),\]

You may easily check the above formulas define a comodule map $\Ad: \gl  \rightarrow \gl \otimes_A k[G]$. We want to calculate the matrix corresponding to $\Ad$ in the basis $B$. The map $\Ad$ corresponds $1-1$ to a map (also denoted $\Ad$)

\[ \Ad: \gl  \otimes_A k[G] \rightarrow \gl  \otimes_A k[G]. \]

For any $C \in \alg$ we let the module $\gl \otimes_A C$ have the basis $B_{C}:= \{ e_{ij} \otimes 1_{C} \}$. Here $1_{C} \in C$ is the multiplicative unit. We get the following matrix:

\[ 
[\Ad]^{B_{k[G]}} _{B_{k[G]}}=
\begin{pmatrix}   \x\w/D  &  -\x\z/D  & \y\w/D &   -\y\z/D \\
                             -\x\y/D  & \x^2/D &  -\y^2/D &  \x\y/D   \\
                             \z\w/D  &  -\z^2/D &  \w^2/D   &  -\z\w/D     \\
                            -\y\z/D  & \x\z/D   &  -\y\w/D  & \x\w/D            
\end{pmatrix}.
\]

Let $I_1:=(\y,\z), I_2:=(\x,\w)$ and $I:=I_1I_2$. Let $k[N]:=k[G]/I, k[T]:=k[G]/I_1$. It follows $T \subseteq G$ is the diagonal torus in $G$ and $T \subseteq N$ is the \emph{schematic normalizer} of $T$.
We get an induced comodule 

\[ \Ad_N : \gl  \otimes k[N] \rightarrow \gl  \otimes k[N]  \]

 defined as follows:

\[ 
[\Ad_N]^{B_{k[N]}} _{B_{k[N]}}=
\begin{pmatrix}    (1,0) & 0 & 0 & (0,1)  \\
                                0 & t_1t_2^{-1} & u_1u_2^{-1} & 0    \\
                              0 & u_1^{-1}u_2 & t_1^{-1}t_2 & 0      \\
                             (0,1) & 0 & 0 & (1,0)             
\end{pmatrix}.
\]

We may pass to the diagonal torus $T$ and get a comodule map 

\[  \Ad_T: \gl \otimes  k[T] \rightarrow \gl   \otimes_A k[T] \]

defined as follows:

\[ 
[\Ad_T]^{B_{k[T]}} _{B_{k[T]}}=
\begin{pmatrix}     1 & 0 & 0 &  0   \\
                                0 & t_1t_2^{-1} &  0  & 0    \\
                              0 &  0  & t_1^{-1}t_2 & 0      \\
                              0  & 0 & 0 & 1             
\end{pmatrix}.
\]

Hence the basis vectors $e_{11}, e_{12}, e_{21}, e_{22}$ have weights $(0,0), (-1,1), (1,-1), (0,0)$. Let $\alpha:=(1,-1)$ and let $R_G:=\{-\alpha, \alpha\}$.

Because of the isomorphim $\Lie(G)^F(-) \cong (\gl)_a(-)$ of functors there is an isomorphism $\Lie(\GL) \cong \gl $.

Let $\Lie(T):=A\{e_{11}, e_{22} \}$ and $\Lie(\GL)_{R_G}:=A\{e_{12}, e_{21} \}$. From the above calculation we get comodul structures

\[ \Ad_{R_G}: \Lie(\GL)_{R_G} \rightarrow \Lie(\GL)_{R_G} \otimes k[N]  \]

and

\[ \Ad_{T}: \Lie(T) \rightarrow \Lie(T) \otimes k[N]  \]

and a direct sum of $k[N]$-comodules $\Lie(\GL) \cong \Lie(T) \oplus \Lie(\GL)_{R_G}$.

\begin{definition}  \label{roots} The direct sum decomposition  

\[ \Lie(\GL) \cong \Lie(T) \oplus \Lie(\GL)_{R_G} \] 

is the \emph{refined weight space decomposition} of $\Lie(\GL)$. The set $R_G$ is the \emph{roots of $\GL$}.
\end{definition}

\begin{example} The adjoint representation for $\Lie(\SL)$. \end{example}

Let $\sl:=A\{e_{12} , e_{11}-e_{22}, e_{21}\} \subseteq \gl $, and let $x:=e_{12}, H:=e_{11}-e_{22}, y:=e_{21}$ and let $B(s):=\{x,H,y\}$. Define for any commutative $A$-algebra $C$ the following: Let $B(s)_C:=\{ x \otimes 1_C, H \otimes 1_C, y \otimes 1_C\}$ be a basis for $\sl  \otimes_A C$. Here $1_C \in C$ is the multiplicative unit. You may check the following:

\[ \Ad(x)= x \otimes (\x^2/D) - H \otimes (\x\z/D) - y \otimes (\z^2/D), \]

\[ \Ad(H)= - x \otimes ((2\x\y)/D)    +  H \otimes ((\x\w + \y\z)/D)   +  y \otimes (2\z\w/D) , \]

and 

\[ \Ad(y)= -  x \otimes (\y^2/D)  +  H \otimes (\y\w/D)   +     y \otimes ( \w^2/D) . \]

Hence $\sl  \subseteq \gl $ is a sub $k[G]$-comodule. Let $S:=\SL:=V(D-1) \subseteq \GL$. we get the induced $k[S]$-comodule

\[ \Ad(S): \sl   \rightarrow \sl  \otimes_A k[S] \]

defined as follows:

\[ \Ad(S)(x)= x \otimes (\x^2) - H \otimes (\x\z) - y \otimes (\z^2), \]

\[ \Ad(S)(H)= - x \otimes (2\x\y)    +  H \otimes (\x\w + \y\z))   +  y \otimes (2\z\w) , \]

and 

\[ \Ad(S)(y)= -  x \otimes \y^2  +  H \otimes \y\w   +     y \otimes \w^2  . \]

A calculation shows the elements $x,H,y$ have weights $\{2,0,-2\}$ as expected. Let $R_S:=\{-2,2\}$.

There is an isomorphism $\Lie(\SL) \cong \sl $. Let $\Lie(T):=A\{H\}$ and let $\Lie(\SL)_{R_S}:=A\{x,y\}$. It follow there are comodule structures

\[ \Delta_1: \Lie(T) \rightarrow \Lie(T) \otimes k[N] \]

and

\[ \Delta_2: \Lie(\SL)_{R_S}  \rightarrow \Lie(\SL)_{R_S} \otimes k[N] \]

And a direct sum decomposition of $k[N]$-comodules:

\[ \Lie(\SL) \cong \Lie(T) \oplus \Lie(\SL)_{R_S}.\]

\begin{definition} The set $R_S:=\{-2,2\}$ is the \emph{roots of $\SL$} and the direct sum decomposition 

\[ \Lie(\SL) \cong \Lie(T) \oplus \Lie(\SL)_{R_S}\]

is the \emph{refined weight space decomposition of $\Lie(\SL)$}.
\end{definition}

\section{The algebra of distributions of $\GL$ and $\SL$}

In this section we include some definitions and results on the algebra of distributions $\Dist_n(G,e)$ of a group scheme $G$ over a commutative ring $k$. We also calculate the module of first order distributions and Lie algebras of $\SL$ and $\GL$.
Let $(k[G], \Delta, S, \epsilon)$ be the Hopf algebra of the group scheme $G$ and let $I:=ker(\epsilon)$ 
be the kernel of the counit - the augmentation ideal. 

\begin{definition} Let $G$ be an affine group scheme over $k$ with augmentation ideal $I \subseteq k[G]$.
We define the $k$-module of distributions of order $n$ as follows: $\Dist_n(G,e):=\Hom_k(k[G]/I^{n+1}, k)$. 
We define $\Dist_n^+(G,e):= \Hom_k(I/I^{n+1},k)$ and $\Lie(G):=\Dist^+_1(G,e):=\Hom_k(I/I^2,k)$.
\end{definition}

For any set of $n$ elements $x_1,x_2,..,x_n$ in $I$ there is the following formula (see \cite{jantzen}, Section I.7.7):

\[\Delta(x_1\cdots x_n)= \prod_{i=1}^n(1 \otimes x_i + x_i \otimes 1) + z \]

with $z \in \sum_{r=1}^n I^r \otimes I^{n+1-r}$. If $n=2$ we get the formula

\[ \Delta(x_1x_2)=(1 \otimes x_1 + x_1 \otimes 1)(1\otimes x_2 + x_2 \otimes 1) + z =\]

\[ 1 \otimes x_1x_2 + x_1 \otimes x_2 + x_2 \otimes x_1 + x_1x_2 \otimes 1 + z \]

with $z \in I \otimes I^2 + I^2 \otimes I \subseteq k[G] \otimes k[G]$.  The counit $\epsilon$ has a section $s: k \rightarrow k[G]$ with $\epsilon \circ s =Id_k$ equal to the identity. Let $\phi:= s \circ \epsilon$. There is an exact sequence

\[ 0 \rightarrow I/I^2 \rightarrow k[G]/I^2 \rightarrow^{\epsilon_1}  k \rightarrow 0 \]

which has a section $s_1$. Let $\phi_1:= s_1 \circ \epsilon_1$.

\begin{lemma} \label{splitting} The section $s_1$ gives a direct sum decomposition 

\[  p_1: k[G]/I^2 \cong I/I^2 \oplus k \]

defined by $p_1(z):=(z-\phi_1(z), \epsilon_1(z))$. The inverse mep $p_2: I/I^2 \oplus k \rightarrow k[G]/I^2$ is defined by $p_2(x,a):=x+s_1(a)$. We get a $k$-linear projection map $q: k[G]/I^2 \rightarrow I/I^2$

\end{lemma}
\begin{proof} The proof follows since $\phi_1 \circ \phi_1 =\phi_1$.
\end{proof}

From Lemma \ref{splitting} we get a sequence

\[ k[G] \rightarrow^{\Delta} k[G] \otimes k[G] \rightarrow k[G]/I^2\otimes k[G]/I^2 \rightarrow I/I^2 \otimes I/I^2 ,\]

where the rightmost map is the map $q \otimes q$. We get a  canonical map

\[ \rho: k[G] \rightarrow k[G]/I^2 \otimes k[G]/I^2.\]

For any pair $u,v \in \Dist_1(G,e)$ we get two maps $u\overline{\otimes} v, v\overline{ \otimes} u: k[G]/I^2 \otimes k[G]/I^2 \rightarrow k \otimes k \cong k$. We may also define the difference

\[ \psi(u,v):= u\overline{\otimes} v- v\overline{ \otimes} u: k[G]/I^2 \otimes k[G]/I^2 \rightarrow k.\]

We get a composed map

\[ k[G ] \rightarrow^{\rho} k[G]/I^2 \otimes k[G]/I^2 \rightarrow^{\psi(u,v)} k.\]

\begin{lemma} \label{welldefined} For any element $z \in I^2 \subseteq k[G]$ it follows $\psi(u,v)(\rho(z))=0$. Hence we get an induced map $\phi(u,v):k[G]/I^2 \rightarrow k$.  It follws $\phi(u,v) \in \Dist_1(G,e)$ is a distribution of order $1$.
For any pair $u,v\in \Dist_1(G,e)$ we define $[u,v]:=\phi(u,v) \in \Dist_1(G,e)$. This gives $\Dist_1(G,e)$ the structure of a $k$-Lie algebra. There is an induced structure as $k$-Lie algebra on $\Lie(G)$.
\end{lemma}
\begin{proof} For any pair of elements $x,y \in I$  we get the following:

\[ \Delta(xy)=1\otimes xy + x\otimes y + y \otimes x + xy \otimes 1 + \]

with $z \in I \otimes I^2 + I^2 \otimes I$.  It follows

\[ \rho(xy)= \overline{x} \otimes \overline{y} + \overline{y} \otimes \overline{x} \in k[G]/I^2 \otimes k[G]/I^2.\]

It follows 

\[ \phi(u,v)(xy)=\psi(u,v)(\rho(xy))=   u(\overline{x})v( \overline{y}) + u(\overline{y})v( \overline{x})  \]

\[   -  v(\overline{x})u(\overline{y}) -v(\overline{y})u( \overline{x})=0.\]

We get an induced $k$-linear map

\[ \phi(u,v):k[G]/I^2 \rightarrow k.\]

If $u,v \in \Lie(G)$ it follows  $\phi(u,v)(1)=0$ hence $\phi(u,v ) \in \Dist^+_1(G,e):=\Lie(G)$. One checks this gives $\Dist_1(G,e)$ and $\Lie(G)$ the structure as a $k$-Lie algebra(s) and the Lemma follows from \cite{jantzen}, CH I.

\end{proof}

Note: In \cite{jantzen} it is proved that the full algebra of distributions $\Dist(G,e)$ is an associative unital $k$-algebra equipped with a filtration $F_n:=\Dist_n(G,e)$ such that the associated graded algebra $Gr(\Dist(G,e),F_n)$ is a commutative unital $k$-algebra.
If $G$ is a group scheme of finite type over a field $k$ of characteristic zero it follows $\Dist(G,e) \cong \U(\Lie(G))$ equals the universal enveloping algebra of the Lie algebra $\Lie(G)$.
If $k$ is a field of characteristic $p>0$ it follows the restricted universal enveloping algebra $\U^{[p]}(\Lie(G))$ embeds into $\Dist(G,e)$ as a strict sub algebra. Hence the algebra of distributions $\Dist(G,e)$ is larger than $\U^{[p]}(\Lie(G))$.
Over field $k$ of characteristic $p>0$ the algebra of distributions plays the role of the universal enveloping algebra. For any integral group scheme $G$ over $k$ it follows for any pair of $k[G]$-comodules $U,V$, there is an "equality"

\[  \Hom_{k[G]-mod}(U,V) = \Hom_{\Dist(G,e)}(U,V), \]

hence we may use the algebra of distributions to make explicit calculations. We may also use the algebra of distributions to study the weight space decomposition of $U,V$.

Note moreover: There is no "equality"

\[  \Hom_{k[G]-mod}(U,V) = \Hom_{\U^{[p]}(\Lie(G))}(U,V), \]

hence you cannot in general use the restricted universal enveloping algebra $\U^{[p]}(\Lie(G))$ to study the full category of $k[G]$-comodules over a field of characteristic $p>0$. 

If  $G$ is a group scheme of finite Krull dimension over $k$ it follows the restricted universal enveloping algebra is finite dimensional over $k$. The algebra of distributions $\Dist_1(G,e)$ is infinite dimensional over $k$. Hence you cannot reduce the study of the category of $k[G]$-comodules to the study of modules on a finite dimensional algebra - there are some properties of $G$ that are "out of reach" or "inaccessible" from the point of view of the theory of finite dimensional algebras.  You need the full algebra of distributions and the construction in \cite{jantzen}, CH I. 

Let $k[G]:=A[x_{ij}][1/D]$ be the localization of $A[x_{ij}]$ at the determinant $D$ and let $\GL:=\Spec(k[G])$. It follows for any commutative $A$-algebra $B$ the group $h_{k[G]}(B)$ is the group of $2\times 2$-matrices with 
coefficients in $B$ and determinant a unit in $B$. There is a canonical surjective map $k[G] \rightarrow k[S]$ with kernel $(D-1)$ hence the group scheme $\SL:=\Spec(k[S])$ is a closed subgroup scheme of $\GL$. We wil calculate the module of distributions 
$\Dist_1(\GL,e)$ and $\Dist_1(\SL,e)$ and the Lie algebras $\Lie(\GL)$ and $\Lie(\SL)$.

The augmentation ideal $I' \subseteq k[G]$ is the ideal $I':=(x_{11}-1, x_{12}, x_{21} , x_{22}-1)$. We may prove there is an isomorphism of free $A$-modules of rank $5$:

\[  k[G]/(I')^2 \cong A\{v_0,v_1,v_2,v_3,v_4\} \]

where $v_0:=1, v_1:=x_{11}-1, v_2:=, x_{12},v_3:= x_{21}$ and $v_4:= x_{22}-1$.  Let $z_i:=v_i^*$ denote the dual basis. It follows there is an isomorphism of free  $A$-modules of rank $5$:

\[ \Dist_1(\GL,e) \cong A\{z_0,z_1,z_2,z_3,z_4\}. \]

We want to calculate the Lie products of the elements $z_i,z_j$ and to classify the module of distributions as $A$-Lie algebra. Let us calculate $\Delta(v_i)$. We get by definition

\[ \Delta(v_0):= v_0 \otimes v_0.\]

We moreover get

\[ \Delta(x_{11}-1):= (x_{11}-1+1) \otimes (x_{11}-1+1) + x_{12} \otimes x_{21} - 1 =\]

\[  (x_{11}-1) \otimes (x_{11}-1) + (x_{11}-1) \otimes 1 + 1 \otimes (x_{11}-1) + 1 -1 + x_{12} \otimes x_{21} =\]

\[ v_1 \otimes v_1 + v_0 \otimes v_1 + v_1 \otimes v_0 + v_2 \otimes v_3.\]

Similarly we get the following calculation:




We define

\[ \rho'(v_1):= v_1 \otimes v_1 + v_2 \otimes v_3 + v_0 \otimes v_1 + v_1 \otimes v_0 ,\]

\[ \rho'(v_2):= v_1 \otimes v_2 + v_2 \otimes v_4 + v_0 \otimes v_2 + v_2 \otimes v_0 ,\]

\[ \rho'(v_3):= v_3 \otimes v_1 + v_4 \otimes v_3 + v_0 \otimes v_3 + v_3 \otimes v_0 ,\]

and

\[ \rho'(v_4):= v_4 \otimes v_4 + v_3 \otimes v_2 + v_0 \otimes v_4 + v_4 \otimes v_0 .\]

There is by Lemma \ref{welldefined} for each $z_i,z_j$ a well defined map 

\[  \psi(z_i,z_j) \circ \rho': k[G]/(I')^2 \rightarrow A,\]

 and we define $[z_i,z_j]:= \psi(z_i,z_j) \circ \rho' \in \Dist_1(\GL,e)$.   Let us calculate $[z_2,z_3]$ as an example:

We get 

\[ \psi(z_2 ,z_3)(\rho'(v_1))=1,  \]

and $\psi(z_2 , z_3)(\rho'(v_i))=0$ for $i=0,2,3,4$.

We also get

\[ \psi(z_3 , z_2)(\rho'(v_4))=1,  \]

and $\psi(z_3 , z_2)(\rho'(v_i))=0$ for $i=0,1,2,3$. From this it follows $[z_2,z_3]=z_1-z_4$. Similarly one calculates $[z_i, z_j]$ for any $i,j$.

Define the following matrices in the ring $M_2(A)$ of $2\times 2$-matrices with coefficients in $A$:

\[
t_1:=  
\begin{pmatrix}   1  &  0 \\
                             0  & 0 
\end{pmatrix},
\]

\[ 
\tilde{x}:=  
\begin{pmatrix}   0  &  1 \\
                             0  & 0 
\end{pmatrix},
\]

\[ 
\tilde{y}:=  
\begin{pmatrix}   0 &  0 \\
                             1  & 0 
\end{pmatrix},
\]

and

\[ 
t_2:=  
\begin{pmatrix}   0 &  0 \\
                             0  & 1 
\end{pmatrix}.
\]

Let $\gl $ denote the $A$-Lie algebra of $2\times 2$-matrices with coefficients in $A$. It follows $\gl  \cong A\{t_1,t_2,\tilde{x},\tilde{y} \}$ is an isomorphism of $A$-Lie algebras. There is an isomorphism $\sl  \cong A\{\tilde{x}, \tilde{H}, \tilde{y}\}$
with $\tilde{H}:=t_1-t_2$.

 There is a surjection $R' \rightarrow R$ with kernel $(D-1)$ giving a closed embedding of group schemes $\SL \subseteq \GL$ . We get a canonical surjective map $p:R'/(I')^2 \rightarrow R/I^2$ and an induced  embedding of modules of distributions 

\[  g:\Dist_1(\SL,e) \subseteq \Dist_1(\GL,e).\]

The $A$-module $R/I^2$ has basis  $R/I^2 \cong A\{ 1,x_{12}, x_{11}-1, x_{21} \} \cong A\{v_0,e_x,e_H,e_y\}$ with $v_0:=1, e_x:=x_{12}, e_H:=x_{11}-1$ and $e_y:=x_{21}$. The dual $\Dist_1(\SL,e)$ has basis

\[ \Dist_1(\SL,e) \cong A\{ w_0, x,H,y\} \]

with $w_0:=v_0^*$.

\begin{theorem} \label{dist} The following holds: $[z_0,z_i]=0$ for $i=1,2,3,4$. Moreover $[z_1,z_2]=z_2, [z_1,z_3]=-z_3, [z_1,z_4]=0, [z_2,z_3]=z_1-z_4, [z_2,z_4]=z_2$ and  $[z_3,z_4]=-z_3$.
Hence there is an isomorphism $f:\Dist_1(\GL,e) \cong A\{z_0\} \oplus \gl $  defined by $f(z_1):=t_1, f(z_2):= \tilde{x}, f(z_3):=\tilde{y}$ and $f(z_4):= t_2$.
The embedding $g$ is defined as follows: $g(w_0):=z_0, g(x):=z_2=\tilde{x}, g(H):=z_1-z_4=\tilde{H}$ and $g(y):=z_3=\tilde{y}$. Hence it identifies $\Dist_1(\SL,e)$ with the sub Lie algebra $A\{z_0\}\oplus \sl $ of  $\Dist_1(\GL,e)$ under the isomorphism $f$.
\end{theorem}
\begin{proof}  The proof is an elementary calculation using the calculations and constructions  above.
\end{proof}

Hence the $A$-Lie algebras $\Lie(\SL)$ and $\Lie(\GL)$ are "what we expect them to be".

\section{On the schematic normalizer of the diagonal torus in $\SL$}

In this section we study the group scheme $G:=\SL$ for any commutative ring $A$ and the schematic normalizer $N \subseteq \SL$ of the diagonal torus. The diagonal torus 
$T \subseteq \SL$ is contained in the schematic normalizer $N$ and the schematic normalizer $N$ is isomorphic to the functor normalizer $N_G(T)$ as defined in the book of Jantzen \cite{jantzen}.


Let $A$ be any commutative unital ring and let $x_{11}, x_{12}, x_{21}, x_{22}$ be independent varaibles over $A$. Let $D:=x_{11}x_{22}-x_{12}x_{21}$ be the "determinant" and let $F(x_{ij}):=D-1 \in A[x_{ij}]$. Let $R:=A[x_{ij}]/(F)$. Let $\alg$ be the category of commutative $A$-algebras and let $\grp$ be the category of groups. Define the functor

\[ h_R: \alg \rightarrow \grp \]

by $h_R(B):=\Hom_{A-alg}(R,B)$. It follows $h_R(B)$ is the group of $2 \times 2$ matrices with coefficients in $B$ and determinant $1$. It follows $G:=\SL:=\Spec(R)$ is an affine group scheme in the sense of the book of Demazure/Gabriel and Jantzen and the ring $R$ is a Hopf algebra. 

Let $I_1:=(x_{12}, x_{21}), I_2:=(x_{11}, x_{22})$ and $I:=I_1I_2$ be ideals in $R$ and let $k[N]:=R/I, k[T]:=R/I_1$.  Let $N:=\Spec(k[N]), T:=\Spec(k[T])$.   Let $t:=\overline{x_{11}}, \frac{1}{t}:=\overline{x_{22}}, u:=\overline{x_{12}}, v:=\overline{x_{21}} \in R/I$. 

\begin{lemma} It follows $R/I \cong R_1\oplus R_2$ with $R_1:=A[t, \frac{1}{t}]$ and $R_2:=A[u,v]/(uv+1)$. The schemes $T,N$ are closed subgroup schemes of $\SL$. The group scheme $T$ is a closed and normal subgroup scheme of $N$.

\end{lemma}
\begin{proof} The comultiplication $\Delta$ satisfies 

\[ \Delta(I) \subseteq R \otimes I + I \otimes R \]

an

\[ \Delta(I_1) \subseteq R \otimes I_1 + I_1 \otimes R \]

and it follows $T,N \subseteq \SL$ are closed subgroup schemes. There is an equality $I_1 + I_2 =(1)$, hence the two ideals $I_1,I_2$ are coprime. It follows $I:=I_1I_2 = I_1 \cap I_2$ annd by the Chinese remaider Lemma
there is an isomorphism of rings

\[  R/I \cong R_1 \oplus R_2, \]

where $R_i:=R/I_i$ for $i=1,2$. The Lemma follows.
\end{proof}

\begin{definition}  The scheme $T$ is the \emph{diagonal torus} in $\SL$ and the scheme $N$ is the \emph{schematic normalizer} of $T$.
\end{definition}


\begin{proposition}  \label{functornorm} There is an isomorphism of group schemes $N \cong N_G(T)$ where $N_G(T)$ is the functor normalizer as defined in \cite{jantzen}.
\end{proposition}
\begin{proof} Let $T \subseteq S:=\SL$ be the diagonal torus. By definition it follows $T:=V(I_1)$ is defined by the ideal $I_1$. There is an equality of ideals $I=(x_{11}x_{12}, x_{21}x_{22})$. By definition it follows $h_{N_G(T)}(R)$ is the group of all elelments $g \in h_S(R)$ with the property that for any map $f: R \rightarrow R_1$ and any element $t \in h_T(R_1)$,  it follows $gtg^{-1} \in h_T(R_1)$. If $g \in V(I_1)$ it follows $g \in h_{N_G(T)}(R)$. Conversely if $g \in h_{N_G(T)}(R)$ the following holds for any element $t \in h_T(R_1)$ with $R_1:=R[z, 1/z]$: $gtg^{-1} \in h_T(R_1)$. Let $g:=(a_{ij})$ be the 2 times 2 matrix with coefficients in $R$ and determinant $D:=a_{11}a_{22}-a_{12}a_{21}$ a unit in $R$ corresponding to $g$.
Let  $t\in h_T(R_1)$ be the element

\[ 
t:=  
\begin{pmatrix}   z &  0 \\
                             0  & z^{-1} 
\end{pmatrix}.
\]

It follows 

\[ 
gtg^{-1}:=  
\begin{pmatrix}   (a_{11}a_{22}z-a_{12}a_{21}/z)\frac{1}{D} &  a_{11}a_{12}(1/z-z)\frac{1}{D}  \\
                             a_{21}a_{22}(z-1/z)\frac{1}{D}  &   (a_{11}a_{22}/z-a_{12}a_{21}z)\frac{1}{D}  
\end{pmatrix}.
\]

If $gtg^{-1} \in h_T(R_1)$ it follows $a_{11}a_{12}=a_{21}a_{22}=0$ and this implies $g \in V(I_1)$. Hence $g\in h_N(A)$ and the Proposition follows.
\end{proof}

The comultiplication $\Delta_N$ on $k[N]$ is the following map:

\[ \Delta_N((t,0)):=  (t,0)\otimes (t,0)+ (0,u)\otimes (0,v), \]

\[ \Delta_N((1/t,0)):=  (0,v)\otimes (0,u)+ (1/t,0 )\otimes (1/t,0), \]

\[ \Delta_N((0,u)):=  (t,0)\otimes (0,u)+ (0,u)\otimes (1/t,0), \]

and 

\[ \Delta_N((0,v)):=  (0,v)\otimes (t,0)+ (1/t,0)\otimes (0,v) \]

We get an induced comodule structure

\[ \Delta_1: k[N] \rightarrow k[N] \otimes_A k[T]  \]

defined by

\[  \Delta_1((t,0)):= (t,0) \otimes t, \]

\[  \Delta_1((1/t,0)):= (1/t,0) \otimes 1/t ,\]

\[  \Delta_1((0,u)):= (0,u) \otimes 1/t ,\]

and

\[  \Delta_1((0,v)):= (0,v) \otimes t .\]

Define the following subring of $k[N]$: 

\[ k[N/T]:=\{ x\in k[N] : \Delta_1(x)=x \otimes 1 \in k[N] \otimes_A k[T] \}, \]

and let $N/T:=\Spec(k[N/T])$. Let $i \in \ZZ$ be an integer. An element $(f(t), g(u,v)) \in k[N]$ is on the following form: 

\[  f(t):=\sum_{i=-m}^m \alpha_i t^i, \]

and since $v:=-1/u$ we may express $g(u,v)$ as a Laurent polynomial in the $u$-variable: We get

\[ g(u,v) =\sum_{j=-m}^m \beta_j u^j.\]

Here $\alpha_i, \beta_j \in A$.

We get the following calculations:

 \[  \Delta_1((t^i,0)):= (t^i,0) \otimes t^i, \]

and

\[  \Delta_1((0 ,u^i)):= (0,u^i) \otimes t^{- i}. \]

If $x:=(f(t), g(u)) \in k[N]$ we get the following:

\[ \Delta_1(x)= \sum_{i=-m}^m \alpha_i(t^i,0)\otimes t^i + \sum_{j=-m}^m \beta_j(0,u^j)\otimes t^{-j} = \]

\[ \sum_{i=-m}^m (\alpha_i t^i, \beta_{-i} u^{-i}) \otimes t^i .\]

And it follows $\Delta_1(x)=x \otimes 1$ iff $x= (\alpha_0, \beta_0)$ with $\alpha_0, \beta_0 \in A$. Hence $k[N/T] \cong A\{e_1,e_2\}$ with $e_1:=(1,0), e_2:=(0,1) \in R_1 \oplus R_2$.  The inclusion $k[N/T] \subseteq k[N]$ induce 
the following structure as Hopf algebra on $k[N/T]$:

There is a comultiplication map 

\[ \Delta_{N/T}: k[N/T] \rightarrow k[N/T] \otimes_A k[N/T] \]

defined by 

\[ \Delta_{N/T}(e_1):= e_1 \otimes e_1 + e_2 \otimes e_2, \]

and

\[ \Delta_{N/T}(e_2):= e_1 \otimes e_2 + e_2 \otimes e_1. \]

The inversion map $S_{N/T} :k[N/T] \rightarrow k[N/T]$ is the identity map. The counit $\epsilon_{n/T}: k[N/T] \rightarrow A$ is the following map:

\[ \epsilon_{N/T}(e_1):= 1,  \epsilon_{N/T}(e_2):= 0. \]

One verifies the inclusion map $k[N/T] \subseteq k[N]$ is a map of Hopf algebras.

Define for any integer  $i \in \ZZ$ the following element: $w_i:=(t^i, u^i) \in k[N]$. 

\begin{lemma}\label{ffmap} The set $B:=\{w_i: i\in Z\}$ is a basis for $k[N]$ as a free module over the ring $k[N/T ]$. It follows the ring extension $k[N/T] \subseteq k[N]$ is faithfully flat.
\end{lemma}
\begin{proof} If $x:=(f(t), g(u)) \in k[N]$ is any element we may write 

\[ f(t):=\sum_{i=-m}^m a_i t^i , g(u):= \sum_{i=-m}^m b_i u_i \]

and it follows

\[ x:=(f(t), g(u))= \sum_{i=-m}^m (a_i,b_i)(t^i, u^i):= \sum_i (a_i,b_i)w_i .\]

One checks the elements $w_i$ is a basis for $k[N]$ as $k[N/T]$-module and the Lemma follows.

\end{proof}

Let $\sigma, q: N \times T \rightarrow N$ be the action morphism and projection morphism.

\begin{proposition}\label{quotient}    The affine scheme $N/T$ is a group scheme and there is a faithfully flat morphism of group schemes $\pi: N \rightarrow N/T$ giving an exact sequence of group schemes 

\[ e \rightarrow T \rightarrow N \rightarrow^{\pi}  N/T \rightarrow e\]
 
with $T \cong ker(\pi)$. The group scheme $N/T$ satisfies the following universal property: For any affine scheme $Z$ and any morphism $\phi:N \rightarrow Z$ with $ \phi \circ \sigma = \phi \circ q$, it follows there is a unique
morphism $ \phi^*: N/T \rightarrow Z$ with $ \phi^* \circ \pi = \phi$.

\end{proposition}
\begin{proof} The inclusion map $i:k[N/T] \rightarrow k[N]$ is a map of Hopf algebras and we get a corresponding map $\pi: N \rightarrow N/T$ of group schemes with $ker(\pi)=T$. It follows from Lemma \ref{ffmap} the map $\pi$ is a surjection since $i$ is faithfully flat. Given any map $\phi: N \rightarrow Z:=\Spec(C)$ with corresponding ring map $f: C \rightarrow k[N]$ 
it follows $\Delta_1(f(x))=f(x) \otimes 1 \in k[N]$ iff the map $f$ factors into a map $f^*: C \rightarrow k[N/T]$ with $f^* \circ i =f$. Hence there is a unique map of affine schemes $\phi^*: N/T \rightarrow Z$ with $\phi^* \circ \pi = \phi$. The Proposition follows.
\end{proof}

By the universal property from Proposition \ref{quotient} it follows $N/T$ is "uniquely determined up to isomorphism". This is why we make the following "definition":

\begin{definition} We say the group scheme $N/T$ is the \emph{quotient of $N$ by the subgroup scheme $T$}.
\end{definition}


Note: If  $G$ is a group scheme of finite type over a Dedekind domain and if $H \subseteq G$ is a closed subgroup scheme, it follows there is a well defined (meaning "unique up to isomorphism") quotient scheme $\pi: G \rightarrow G/H$ satisfying the universal property from Proposition \ref{quotient}.
If $G$ is a group scheme over a more general base $U$, it follows the quotient $G/H$ may no longer be a scheme. In such situations one has access to the notions \emph{algebraic space} and \emph{algebraic stack}.
In \cite{DG}, Theoreme III.2.3.2 there is a general criterium for a quotient of a groupoid scheme $G_1 /G_0$ to be a scheme. One needs an open affine cover of $G_0$ satisfying a compatibility property.  In the case when $G$ is a finite group scheme acting on a quasi projective scheme $X \subseteq \mathbb{P}^n_k$ of finite type over a field $k$, it follows there is always an open affine $G$-invariant cover $U_i$ of $X$. In this case we may use the open affine cover $\{U_i\}$  to define the 
quotient scheme $X/G$. 

The reason for the inclusion of Proposition \ref{quotient}
is the explicit construction of the Hopf algebra $k[N/T]$ and the explicit construction of the quotient map $\pi$. This construction will be needed in the paper. In \cite{waterhouse}, Lemma 16.6.3 you will find an explicit and elementary construction of the quotient $G/N$ for any affine group scheme $G$ of finite type over a field and any closed normal subgroup scheme $N \subseteq G$. The above Proposition \ref{quotient} holds over any commutative base ring $A$. The proof of Proposition \ref{quotient} differs from the proof of the Lemma in \cite{waterhouse} - the proof in \cite{waterhouse} needs the base ring to be a field.

\begin{example} Finite group schemes, semi direct products and the Weyl group. \end{example}


 The exact sequence

\[ e \rightarrow T \rightarrow N \rightarrow N/T \rightarrow e \]

does not split in general, hence the group scheme $N$ is not the semi direct product of $T$ with $N/T$.  One may check the following explicitly: Let $\C$ be the field of complex numbers and let $\SLO(\C^2)$ be the linear algebraic group of $2 \times 2$-matrices with complex coefficients. Let  $T_{\C} \subseteq \SLO(\C^2)$ be the diagonal torus. We get an exact sequence of abstract groups

\[ e \rightarrow T_{\C} \rightarrow N_{\C} \rightarrow N_{\C}/T_{\C} \rightarrow e \]

and one checks $N_{\C}/T_{\C} \cong S(2) $ is the symmetric group on two letters. One may calculate the cohomology class $c(N_{\C})  \in \H^2(S(2), T)$ associated to the above exact sequence and one will find the class $c(N_{\C})$ is non trivial.
Hence the sequence does not split. The Weyl group $W$ of $\SLO(2,\C)$ is isomorphic to $S(2)$.

\section{On the schematic normalizer in $\GL$}

In this section we study the group scheme $G:=\GL$ where $A$ is an arbitrary commutative unital ring, the diagonal torus $T \subseteq G$ and the schematic normalizer $N \subseteq G$ of $T$. We prove there is a split exact sequence of affine group schemes

\[ e \rightarrow T \rightarrow N \rightarrow^{\pi}  N/T \rightarrow e \]

where $\pi$ i a faithfully flat morphism and $T \cong ker(\pi)$. The affine group scheme $N/T$ satisfies the following universal property: Let $\sigma, q: N \times T \rightarrow N$ be the action map and projection map. For any morphism $f: N \rightarrow Z$ to an affine scheme $Z$ 
with $ \sigma \circ f = q \circ f$, there is a unique morphism $\tilde{f}: N/T \rightarrow Z$ with $\tilde{f} \circ \pi = f$. Hence the affine group scheme $N/T$ is uniquely determined up to isomorphism.
There is an isomorphism of group schemes

\[  N \cong T \times^{\sigma} N/T, \]

where the right hand "product" is the semi direct product with respect to a section $\sigma$ of $\pi$.

Let $R:=A[x_{ij}, 1/D]$ where $D:=x_{11}x_{22}-x_{12}x_{21}$ is the determinant and let $G:=\GL:=\Spec(R)$. Let $I_1:=(x_{12}, x_{21}), I_2:=(x_{11}, x_{22})$ and $I:=I_1I_2$. It follows $I_1 + I_2 =(1)$ and we get an isomorphism
of rings $R/I \cong R/I_1 \oplus R/I_2:=R_1 \oplus R_2$ with 

\[ R_1 \cong A[t_1, 1/t_1, t_2, 1/t_2], R_2:= A[u_1,1/u_1, u_2, 1/u_2] \]

where $\overline{x_{11}}:=(t_1,0) , \overline{x_{22}}:=(t_2,0)   ,   \overline{x_{12}}:=(0 ,u_1) ,  \overline{x_{21}}:=(0 ,u_2) $, and  where $t_i, u_j$ are independent variables over $A$.
Let $N:=V(I) \subseteq G$ and let $T:=V(I_1) \subseteq N$.

\begin{lemma} The group scheme $T$ is a normal subgroup scheme of $N$. There is an isomorphism of group schemes $N \cong N_G(T)$ where $N_G(T)$ is the functor normliazer as defined in \cite{jantzen}.
\end{lemma}
\begin{proof} The first claim is a straightforward calculation. The second claim is similar to Proposition \ref{functornorm}.
\end{proof}

Let for any pair $(i,j) \in \ZZ \oplus \ZZ$  $w(i,j):=(t_1^it_2^j, u_1^ju_2^i)$ and $z(i,j):=(t_1^it_2^j, u_1^iu_2^j)$. Let $e_1:=(1,0), e_2:=(0,1) \in R_1 \oplus R_2$ and let $i: A\{e_1,e_2\} \rightarrow R/I$ be the canonical inclusion map.

\begin{lemma} \label{basisGL} The set $B:=\{w(i,j): (i,j) \in \ZZ \oplus \ZZ\}$ is a basis for $R/I$ as a right module on $A\{e_1,e_2\}$. Hence the canonical inclusion map $i: A\{e_1,e_2\} \rightarrow R/I$ is faithfully flat.
\end{lemma}
\begin{proof} Let $x:=(f, g) \in R/I$ be any element. We may write 

\[   f(t_1,t_2)= \sum_{(i,j) \in S_1} \alpha_{i,j}t_1^it_2^j  \]

and

\[  g(u_1,u_2):= \sum_{(k,l)\in S_2} \beta_{k,l}u_1^k u_2^l \]

for finite subsets $S_1,S_2 \subseteq \ZZ^2$. We may write 

\[  x:=(f,g)=(f,0)+(0,g) =\]

\[x=   \sum_{S_1} (\alpha_{i,j},0)(t_1^it_2^j,0) + \sum_{S_2} (0, \beta_{k,l})(0 , u_1^ku_2^l) \]

and we get

\[   x = \sum_{S_1} (\alpha_{i,j},0)(t_1^i t_2^j , u_1^ju_2^i ) + \sum_{S_2} (0, \beta_{k,l})(t_1^lt_2^k, u_1^ku_2^l). \]

Let $a_{i,j}:=(\alpha_{i,j},0), b_{l,k}:=(0, \beta_{k,l}).$ It follows  $a_{i,j}, b_{l,k} \in A\{e_1,e_2\}$ and 

\[ x= \sum_{S_1} a_{i,j} w(i,j)  + \sum_{S_2} b_{l,k} w(l,k), \] 

and it follows $x$ is a linear combination of vectors in $B$ with coefficients in $A\{e_1,e_2\}$. Hence the set $B$ generates $R/I$ as right module on $A\{e_1,e_2\}$.
 One checks the set $B$ is linearly independent and the Lemma follows.
\end{proof}

We get an induced comodule map

\[ \Delta_N: k[N] \rightarrow k[N] \otimes k[N] \]

defined as follows:

\[ \Delta_N(t_1,0):=\Delta_N(x_{11})= (t_1,0) \otimes (t_1,0) + (0,u_1) \otimes (0,u_2), \]

\[ \Delta_N(t_2,0):=\Delta_N(x_{22})= (t_2,0) \otimes (t_2,0) + (0,u_2) \otimes (0,u_1),  \]

\[ \Delta_N(0,u_1):=\Delta_N(x_{12})= (t_1,0) \otimes (0 ,u_1) + (0,u_1) \otimes (t_2, 0),  \]

and

\[ \Delta_N(0 ,u_2):=\Delta_N(x_{21})= (t_2,0) \otimes (0 ,u_2) + (0,u_2) \otimes (t_1, 0). \]

In general we get the following:

\[ \Delta_N(t_1^kt_2^l , 0 ) = (t_1^k t_2^l ,0) \otimes (t_1^k t_2^l  , 0 ) + (0,u_1^k u_2^l ) \otimes ( 0 ,  u_1^l u_2^k ), \]

and

\[ \Delta_N(0, u_1^ku_2^l ) = (t_1^k t_2^l ,0) \otimes (0 , u_1^k u_2^l ) + (0,u_1^k u_2^l ) \otimes ( t_1^l t_2^k , 0 ). \]

When we pass to the quotient $k[T]$ we get the following comodule map:

\[ \sigma: k[N] \rightarrow k[ N] \otimes k[T]  \]

defined by 

\[  \sigma_T(t_1,0):= (t_1,0) \otimes t_1, \]

\[  \sigma_T(t_2,0):= (t_2,0) \otimes t_2, \]

\[  \sigma_T( 0 ,u_1):= (0 ,u_1) \otimes t_2, \]

and

\[  \sigma_T(0 ,u_2):= (0, u_2) \otimes t_1. \]

We moreover get

\[  \sigma_T(t_1^k t_2^l ,0):= (t_1^kt_2^l ,0) \otimes t_1^kt_2^l , \]

and

\[  \sigma_T(0, u_1^ku_2^l ):= (0, u_1^ku_2^l ) \otimes t_1^lt_2^k. \]

It follows

\[ \sigma_T(w(i,j))= \sigma_T(t_1^it_2^j,0)+ \sigma_T(0, u_1^ju_2^i):=\]

\[ (t_1^it_2^j,0)\otimes t_1^it_2^j +(0, u_1^ju_2^i)\otimes t_1^it_2^j = w(i,j) \otimes t_1^i t_2^j .\]

Since $k[T] \cong_{(k,l) \in \ZZ^2}  A\{t_1^k t_2^l \}$ is a free $A$-module we get an isomorphism

\[  k[N] \otimes_A k[T] \cong \oplus_{(k,l) \in \ZZ^2} k[N]\{t_1^kt_2^l\}. \]

Define the following subring of $k[N]$: 

\[  k[N/T]:=\{x \in k[N]: \sigma(x)=x \otimes 1 \in k[N]\otimes k[T] \}. \]

We may write for any element $x \in k[N]$ the following:

\[ x= \sum_{(i,j)\in \ZZ^2} a_{i,j}w(i,j) \]

with $a_{i,j} \in A\{e_1,e_2\}$, since $B$ is a basis for $k[N]$. We get

\[ \sigma_T(x)= \sum_{(i,j)} a_{i,j}w(i,j) \otimes t_1^it_2^j \]

and it follows $\sigma_T(x)=x \otimes 1$ iff $x=uw(0,0):=a_1e_1+ a_2e_2 \in A\{e_1, e_2\}$. Hence there is an "equality"  $k[N/T] = A\{e_1,e_2\}$.

\begin{lemma} There is an equality $k[N/T]=A\{e_1,e_2\}$. The Hopf algebra structure on $k[N]$ induce the following structure on $k[N/T]$: There is a comultiplication $\Delta_{N/T}$, coinversion $S_{N/T}$ and counit $\epsilon_{N/T}$ defined as follows:

\[ \Delta_{N/T}(e_1):= e_1 \otimes e_1 + e_2 \otimes e_2, \Delta_{N/T}(e_2):= e_1 \otimes e_2 + e_2 \otimes e_1 ,\]

\[ S_{N/T}(e_1)=e_1, S_{N/T}(e_2)=e_2 \]

and

\[ \epsilon_{N/T}(e_1)=1, \epsilon_{N/T}(e_2)=0.\]
There is an ismorphism of $A$-algebras $\phi: A[z]/(z^2-z) \cong k[N/T]$ defined by $\phi(a):=(a,a)$ for $a\in A$ and $\phi(z):=e_1:=(1,0), \phi(1-z)=e_2:=(0,1)$.

\end{lemma}
\begin{proof} The first claim follows from the above discussion. Let us calculate the structure as Hopf algebra on $k[N/T]$. There is an isomorphism $\rho: R/I \cong R_1 \oplus R_2$ and the element $z:=x_{11}x_{22}/D$ is mapped to $\rho(z):=e_1 \in k[N/T]$. The element $w:=-x_{12}x_{21}/D$ is mapped to $\rho(w)=e_2$. We check the following: $\Delta_N(D)= D\otimes D \in k[N]$. It follows

\[ \Delta_N(e_1)=  (1/ D)\otimes (1/ D)( x_{11}\otimes x_{11} + x_{12} \otimes x_{21})(x_{21} \otimes x_{12} + x_{22} \otimes x_{22})=\]

\[  e_1 \otimes e_1 + (-e_2) \otimes (-e_2) = e_1 \otimes e_1 + e_2 \otimes e_2.\]

We get simlarly $\Delta_N(e_2)=e_1 \otimes e_2 + e_2 \otimes e_1$. The rest of the Lemma follows with similar calculations.
\end{proof}

The map $\rho: k[N/T] \rightarrow k[N]$ is defined as follows:

\[ \rho(z):= \frac{x_{11}x_{22}}{D} \]

where $D:=x_{11}x_{22} - x_{12}x_{21}$ is the determinant. We get for any commutative $A$-algebra $B$ an induced map

\[ \pi_B: h_N(B) \rightarrow h_{N/T}(B) \]

defined by $\pi_B(g):=\frac{a_{11}a_{22}}{D}$ where  $g\in h_N(B)$ is the following matrix:

\[ 
 g:=  
\begin{pmatrix}   a_{11}   &    a_{12}   \\
                               a_{21}   &   a_{22} 
\end{pmatrix}.
\]

There is a functorial section $\sigma_B: h_{N/T}(B) \rightarrow h_N(B)$ defined as follows: Let $b\in h_{N/T}(B)$ be an element with $b^2=b$ and define

\[ 
 \sigma_B(b):=  
\begin{pmatrix}   b   &    1-b   \\
                               1-b  &   b 
\end{pmatrix}.
\]

It follows $\sigma_B(b_1\times b_1)=\sigma_B(b_1)\sigma_B(b_2)$ and $\sigma_B(1)=Id$. Moreover $ker(\pi_B)=h_T(B)$. Hence the exact sequence

\[  (e) \rightarrow h_T(B) \rightarrow h_N(B) \rightarrow^{\pi_B} h_{N/T}(B) \rightarrow (e) \]

is split by $\sigma_B$. It follows there is an isomorphism

\[  h_N(B) \cong h_T(B) \times^{\sigma_B} h_{N/T}(B) \]

where the right hand side is the semi direct product with respect to $\sigma_B$. The map $\sigma_B$ is functorial and we get a natural transformation 
$\sigma: h_{N/T}(-) \rightarrow h_{N}(-)$ of functors.

Let $\sigma, q: N \times T \rightarrow N$ be the action map and projection map and let $N/T:=\Spec(k[N/T])$.

\begin{theorem} \label{quotientGL} Let $G:=\GL$ with $A$ a commutative unital ring and let $T$ be the diagonal torus with $N \subseteq G$ the schematic normalizer. There  is a faithfully flat map of group schemes $\pi: N \rightarrow N/T$ with $T \cong ker(\pi)$ giving a split exact sequence of group schemes

\[ e \rightarrow T \rightarrow N \rightarrow^{\pi} N/T \rightarrow e .\]

For any map $f: N \rightarrow Z$ to an affine scheme $Z$ with $ f \circ \sigma = f \circ q$ there is a unique map $f^*: N/T \rightarrow Z$ with $f^* \circ \pi = f$. Hence the affine group scheme $N/T$ is uniquely determined up to isomorphism.
The section $\sigma: h_{N/T}(-) \rightarrow h_N(-)$ defined above gives a splitting of the sequence, hence
there is an isomorphism of group schemes $N \cong T \times^{\sigma} N/T$ where the product is the semi direct product with respect to $\sigma$.
 \end{theorem}
\begin{proof} By Lemma \ref{basisGL} it follows the inclusion map $i: k[N/T] \rightarrow k[N]$ is faithfully flat, hence the induced map $\pi: N \rightarrow N/T$ is faithfully flat. It follows $\pi$ is a surjective map of affine group schemes.
One verifies $T \cong ker (\pi)$ equals the kernel of $\pi$. Let $f: C \rightarrow k[N]$ be a map of commutative $A$-algebras with 

\[ f(x):= \sum_{i,j} a_{i,j}w(i,j) \]

for $a_{i,j} \in A\{e_1,e_2\} \cong k[N/T]$. Let $q^*: k[N] \rightarrow k[N] \otimes k[T]$ be  defined as $q^*(c);=c \otimes 1$. It follows  $\sigma_T(f(x))=q^*(f(x)) $ iff 

\[  \sum_{i,j}a_{i,j}w(i,j) \otimes t_1^i t_2^j = \sum_{i,j}a_{i,j}w(i,j) \otimes 1 \]

and this is iff $f(x):= u:= a_1e_1+ a_2e_2 \in k[N/T]$. Hence $f \circ \sigma_T = f \circ q^*$ iff there is a map $f^*: C \rightarrow k[N/T]$ with $f^* \circ i = f$. 
One checks the section $\sigma_B$ defined above for any $B$ gives rise to a natural transformation of functors $\sigma: h_{N/T}(-) \rightarrow h_{N}(-)$ which is a section of the canonical map $\pi:h_N(-) \rightarrow h_{N/T}(-)$.
The Theorem follows.
\end{proof}

Because of Theorem \ref{quotientGL} we say the the affine group scheme $N/T$ is the \emph{quotient of $N$ by $T$}. In \cite{waterhouse} this type of quotient is constructed when $A:=k$ is a field using "other methods".
You may also find a construction on \cite{jantzen} and \cite{DG}.  The reason for the construction given above is the explicit determination of the hopf algebra $k[N/T]$.

\begin{example} The functor normalizer and the Weyl group. \end{example}

\begin{lemma} Let $k[S(2)]$ be the Hopf algebra of the constant group scheme (over $A$) of the symmetric group $S(2)$ on two elements. There is an isomorphism of Hopf algebras $k[S(2)] \cong A[z]/(z(1-z))$, where $A[z]/(z(1-z))$ has the following structure: The comultiplication $\Delta$, coinversion $S$ and conunit $\epsilon$ is defined as follows:   

\[ \Delta(z):= 1 + 2(z\otimes z)-1\otimes z - z \otimes 1, S(z):=1, \epsilon(z):=1. \] 

Hence the group scheme $N/T$ is the constant group scheme of the Weyl group $W_G:=S(2)$ of $\GL$.
\end{lemma}
\begin{proof} The proof is straight forward.
\end{proof}

Hence the functor normalizer $N_G(T) \cong T \times^{\sigma} \Spec(k[W_G])$ is the semi direct product of the diagonal torus $T \subseteq \GL$ and the constant group scheme $\Spec(k[W_G])$ of the Weyl group $W_G$.



\section{A refined Weyl character formula for $\GL$.}


In this section we define and calculate the Weyl group of the group schemes $\GL, \SL$ with $A$ an arbitrary integral domain. We also introduce the \emph{refined weight space decomposition} of the symmetric powers $\Sym^d(V)$
and prove in Theorem \ref{refined} the \emph{refined Weyl character formula}. This formula implies the classical Weyl character formula.

Let $G:=\GL, S:=\SL$ with $A$ an integral domain. Let $T_G \subseteq N_G \subseteq G$ be the diagonal torus and schematic normalizer in $G$ and let 
$T_S \subseteq N_S \subseteq S$ be the diagonal torus and schematic normalizer in $S$. It follows $T_G  \cong \mathbb{G}_m \times \mathbb{G}_m$ and $T_S \cong \mathbb{G}_m$. By definition 

\[ X_G(T):=\Hom(T_G, \gm), Y_G(T):=\Hom(\gm, T_G) \]

and

\[  X_S(T):=\Hom(T_S, \gm), Y_S(T):=\Hom(\gm, T_S).\]

It follows 

\[  X_G(T) \cong \gm \times \gm, Y_G(T) \cong \gm \times \gm, X_S(T) \cong \gm, Y_S(T) \cong \gm.\]

Let $\alpha:=(1,-1) \in X_G(T), \beta:=(-1,1):=-\alpha \in X_G(T)$. We define the \emph{coroots} $\alpha^*, \beta^*$ as follows: Define $I_{\alpha}:=(x_{11}-1, x_{21}, x_{22}-1), I_{\beta}:=(x_{11}-1, x_{12} , x_{22}-1) \subseteq k[G]$. Let $U_{\alpha}:=V(I_{\alpha}), U_{\beta}:=V(I_{\beta}) \subseteq G$. There are \emph{root morphisms}

\[ x_{\alpha}: \ga \cong U_{\alpha},  x_{\beta}: \ga \cong U_{\beta} \]

defined as follows: For any $A$-algebra $B$ define the following:

\[ x_{\alpha}(B): \ga(B) \rightarrow G(B) \]

by

\[ 
x_{\alpha}(B)(b):= 
\begin{pmatrix}   1  &  b \\
                             0  & 1 
\end{pmatrix}.
\]

It follows $x_{\alpha}$ induce an isomorphism $x_{\alpha}: \ga \cong U_{\alpha}$. Define

\[ x_{\beta}(B): \ga(B) \rightarrow G(B) \]

by

\[ 
x_{\beta}(B)(b):= 
\begin{pmatrix}   1  &  0 \\
                             b  & 1 
\end{pmatrix}.
\]

It follows $x_{\beta}$ give an isomorphism $x_{\beta}: \ga \cong U_{\beta}$.

The set $R_G:=\{\beta, \alpha\}$ are the roots of $G$ and there is an isomorphism 

\[ \Lie(G) \cong \Lie(T_G) \oplus (\bigoplus_{\gamma \in R_G} \Lie(G)_{\gamma}). \]

The root subgroups $U_{\gamma}$ are uniquely characterized (up to the choice of a unit in $A$ by \cite{jantzen}, Chapter II) by the following properties: For any $t\in T_G(B)$ and $b \in B$ it follows

\[  tx_{\gamma}(b)t^{-1} = x_{\gamma}(\gamma(t)b) \]

and the induced map

\[  dx_{\gamma}: \Lie(\ga) \cong \Lie(G)_{\gamma} \]

is an isomorphism. For any unit $u \in B^*$ , define the following:

\[  n_{\gamma}(B)(u):= x_{\gamma}(B)(u)x_{-\gamma}(B)(-1/u)x_{\gamma}(B)(u). \]


Define the \emph{coroot} $\gamma^*$ as follows $\gamma^*(B)(u):= n_{\gamma}(B)(u)n_{\gamma}(B)(1)^{-1}.$ It follows

\[  
\alpha^*(B)(u):=
\begin{pmatrix}   u  &  0 \\
                             0  & 1/u 
\end{pmatrix}.
\]

and

\[  
\beta^*(B)(u):=
\begin{pmatrix}  1/ u  &  0 \\
                             0  & u 
\end{pmatrix}.
\]

We get coroots $\alpha^*, \beta^*: \gm \rightarrow T_G$.

There is a canonical perfect pairing 

\[ <,> : X_G(T) \times Y_G(T) \rightarrow \ZZ \]

defined by $<\rho_{m,n}, \gamma_{u,v}>:=mu+nv$, and we define for any root $\gamma \in R_G$ and any element $\rho_{m,n} \in X_G(T)$ the following:

\[  s_{\gamma}(\rho_{m,n}):= \rho_{m,n}-< \rho_{m,n}, \gamma^*> \gamma. \]



\begin{definition} The Weyl group $W_G$ of $G$ is the group of automorphisms of $X_G(T)_{\mathbf{R}}$ generated by $s_{\gamma}$ for $\gamma \in R_G$ a root, and the identity automorphism $e$.
\end{definition}

\begin{lemma}  Let $G:=\GL$.There is an isomorphism $W_G \cong S(2)$ where $S(2)$ is the symmetric group on two elements. The same result holds for $\SL$.
\end{lemma}
\begin{proof} From  Definition \ref{roots} it follows the roots of $G$ are the following:  $\alpha:=(1,-1), \beta=-\alpha=(-1,1)$. By definition we get for any $\rho_{m,n} \in X_G(T)$ the following calculation:

\[  s_{\alpha}(\rho_{m,n}):=\rho_{m,n}-<\rho_{m,n}, \alpha^*> \alpha =\]

\[ \rho_{m,n}-(m-n)\rho_{1,-1} = \rho_{m,n}-\rho_{m-n,n-m}=\rho_{n,m}. \]

Similarly we get $s_{\beta}(\rho_{m,n})=\rho_{n,m}$. Hence the Weyl group $W_G \cong S(2)$ is the symmetric group on two elements. A similar calculation gives the same result for $\SL$.
\end{proof}


Let $(V,\Delta)$ be the standard right comodule on $\GL$. Recall the following: $V:=A\{e_1,e_2\}$ and

\[ \Delta(e_1):= x_{11} \otimes e_1 + x_{12} \otimes e_2, \Delta(e_2):= x_{21} \otimes e_1 + x_{22} \otimes e_2.\]



Let $k[N]:=k[G]/I_1I_2$. It follows $k[N]\cong A[t_1, 1/t_1,t_2,1/t_2]\oplus A[u_1,1/u_1, u_2, 1/u_2]$ with $\overline{x_{11}}:=(t_1,0) , \overline{x_{22}}:=(t_2,0),  \overline{x_{12}}:=(0,u_1)$ and $ \overline{x_{21}}:=(0,u_2)$.

Let $d \geq 1$ be an integer and construct the symmetric power $\Sym^d(V)$. There is a comodule map $\Delta_d: \Sym^d(V) \rightarrow k[G] \otimes \Sym^d(V)$ defined by

\[ \Delta_d(e_1^k e_2^l):= (x_{11} \otimes e_1 + x_{12} \otimes e_2)^k(x_{21} \otimes e_1 + x_{22} \otimes e_2)^l\]

for any pair of integers $0 \leq k,l$ with $k+l=d$. This calculation takes place inside $k[G] \otimes_A \Sym^*(V)$. When we mutilply out we get an element in $k[G]\otimes_A \Sym^d(V)$. Since $(a,0)(0,b)=0$ in $k[N]$ we get the following formula for the induced $k[N]$-comodule

\[ \Delta_{N,d}: \Sym^d(V) \rightarrow k[N]\otimes \Sym^d(V) .\]

We get

\[ \Delta_{N,d}(e_1^ke_2^l):= (x_{11} \otimes e_1 + x_{12} \otimes e_2)^k(x_{21} \otimes e_1 + x_{22} \otimes e_2)^l = \]

\[ ((t_1,0) \otimes e_1 + (0,u_1) \otimes e_2)^k( (0,u_2) \otimes e_1 +  (t_2,0)  \otimes e_2)^l = \]

\[  (t_1^kt_2^l,0) \otimes e_1^ke_2^l + (0, u_1^ku_2^l) \otimes e_1^l e_2^k.\]






We get an induced comodule 

\[ \Delta_{T,d} : \Sym^d(V) \rightarrow k[T] \otimes \Sym^d(V) \]

defined by 

\[ \Delta_{T,d}(e_2^ie_1^{d-i}):= t_1^{d-i}t_2^i \otimes e_2^i e_1^{d-i}.\]

Hence the vector $e_2^ie_1^{d-i}$ has weight $t_2^it_1^{d-i}$ for $i=0,..,d$. Let 

\[  \Gamma(d):=\{ t_2^i t_1^{d-i}: i=0,..,d\} \]

be the weights of $\Sym^d(V)$.

 The Weyl group $W_G \cong S(2)$ acts on $\Gamma(d)$ and the orbits are the sets 

\[  \Gamma(d,i):=\{ t_2^it_1^{d-i}, t_2^{d-i}t_1^i\} \]

 for $i=0,.., k$ with $d:= 2k+1 $. If  $d:=2k \geq 2$ we get the following:

\[  \Gamma(d,i):=\{ t_2^i t_1^{d-i }, t_2^{d-i}t_1^i\} \]

for $i=0,..,k-1$ and $\Gamma(d,k):=\{t_2^k t_1^k\}$. It follows  

\[ \Gamma(d) = \cup_{i=0}^k \Gamma(d,i) \]

 is a disjoint union of sets. It partitions the weights $\Gamma(d)$ into the disjoint union of the orbits under the action of $W_G$.

\begin{definition}  Let $\Sym^d(V)^i:= \oplus_{\lambda \in \Gamma(d,i)} \Sym^d(V)_{\lambda}$. Let $\oplus_{i=0}^k \Sym^d(V)^i$ be the \emph{refined weight space decomposition} of $\Sym^d(V)$.
\end{definition}

By definition there is an isomorphism

\[  \Sym^d(V) \cong \oplus_{i=0}^k \Sym^d(V)^i \]

of $A$-modules.  


We get for any integer $0 \leq i \leq d$  the following calculation in $k[N] \otimes \Sym^d(V)$:

\[  \Delta_{N,d}(e_2^ie_1^{d-i}):=((t_1,0)\otimes e_1 + (0,u_1) \otimes e_2)^{d-i}((0,u_2) \otimes e_1 + (t_2,0)\otimes e_2)^i =\]

\[ (t_1^{d-i}t_2^i,0) \otimes e_1^{d-i}e_2^i + (0, u_1^{d-i}u_2^i) \otimes e_1^i e_2^{d-i}.\]

We get for each $i=0,..,d$ an induced comodule structure

\[  \Delta_{N,d,i} :  \Sym^d(V)^i \rightarrow k[T] \otimes \Sym^d(V)^i \]

and a direct sum of $k[N]$-comodules 

\[ \Sym^d(V) \cong \oplus_{i=0}^d \Sym^d(V)^i.\]

\begin{theorem} \label{refined}(The refined Weyl character formula)  There is  for each integer $i=0,..,d$ a comodule structure

\[ \Delta_{N,d,i}: \Sym^d(V)^{i} \rightarrow k[N]\otimes \Sym^d(V)^{i} \]

and a direct sum decomposition $\Sym^d(V) \cong \oplus_{i=0}^d \Sym^d(V)^{i}$ of $k[N]$-comodules.

The character of $\Sym^d(V)^i$ is the following: If $d:=2k+1 \geq 1$ is an odd integer it follows

\[ Char(\Sym^d(V)^i)=x_2^{d-i}x_1^i + x_2^i x_1^{d-i}\]

For $i=0,..,k$. If $d:=2k \geq 2$  the following holds:

\[ Char(\Sym^d(V)^i)= x_1^{d-i}x_2^i +x_1^i x_2^{d-i} \]

if $i=0,..,k-1$, and 

\[  Char(\Sym^d(V)^{k})=x_1^kx_2^k.\]



\end{theorem}
\begin{proof} The character formula follows from an explicit calculation. The rest of the Theorem follows from the calculation above.
\end{proof}

\begin{corollary}  (The "Weyl character formula"). There is for every odd integer $d  \geq 1$ an equality 

\[ Char(\Sym^d(V))= \sum_{i=0}^k x_2^i x_1^{d-i} + x_2^{d-i} x_1^i. \]
There is for every even integer $d  \geq 2$ an equality

\[ Char(\Sym^d(V)) =  x_1^k x_2^k +\sum_{i=0}^{k-1} x_1^{d-i}x_2^i  + x_1^i x_2^{d-i} .\]

\end{corollary}
\begin{proof} The Corollary follows from Theorem \ref{refined}.
\end{proof}

\begin{example} \label{finitefield} A relation with irreducible comodules over a field of characteristic $p>0$. \end{example}

Let $k$ be a field of characteristic $2$ and let $G:=\SLO_{2,k}$. Let $V \cong k\{e_1,e_2\}$ be the standard right comodule on $G$ and let $\Sym^2(V)$ be the second symmetric power of $V$.
It follows the refined weight space decomposition of $\Sym^2(V)$ is as follows: Define $\Sym^2(V)^1:=k\{e_1^2, e_2^2\}$ and $\Sym^2(V)^2:=k\{e_1e_2\}$. There is a direct sum decomposition of
$k[N]$-comodules

\[ \Sym^2(V) \cong \Sym^2(V)^1 \oplus \Sym^2(V)^2, \]

and the irreducible comodule $V(2)$ with highest weight $2$ is the strict sub-comodule $V(2):=k\{e_1^2,e_2^2\} \subsetneq \Sym^2(V)$. Hence there is an equality 

\[  V(2) = \Sym^2(V)^1 \]

where $\Sym^2(V)^1$ is a term in the refined weight space decomposition. Hence the refined weight space decomposition contains information on irreducible comodules in positive characteristic.

\begin{example} Irreducibility over the complex numbers. \end{example}

Let $k$ be the complex number field and let $G:=\GLO_{2,k}$. Let $(V, \Delta)$ be the standard right comodule and let $\Sym^d(V)$ be the $d$'th symmetric power with $\Sym^d(V)^j$ the $j$'th component in the  refined weight 
space decomposition. The comodule

\[  \Sym^d(V)^j \rightarrow \Sym^d(V)^j \otimes k[N] \]

is an irreducible comodule for all $j$ and $d \geq 1$.








\begin{example} Schur-Weyl functors and refined weight space decompositions.\end{example}

For any Schur-Weyl functor $S_{\lambda}(-)$ as introduced in \cite{fulton} and for the standard right comodule $(V, \Delta)$ on $G:=\GLO_{n,\ZZ}$ we may construct the comodule

\[  \Delta_{\lambda}: S_{\lambda}(V) \rightarrow S_{\lambda}(V) \otimes k[G] \]

and when the Weyl group $W_G$ leaves the weights $\Gamma(\lambda) \subseteq k[T]^*$  of $S_{\lambda}(V)$ invariant we may ask for a map of  $k[N]$-comodules

\[  \Delta_{\lambda, i}:  S_{\lambda}(V)^i  \rightarrow S_{\lambda}(V)^i  \otimes k[N] .\]

The orbits $\Gamma(W,i)$ of the action of $W_G$ on $\Gamma(W)$ defines the direct summand $S_{\lambda}(V)^i:= \oplus_{\gamma \in \Gamma(W,i)} S_{\lambda}(V)_{\gamma}$,  and the  refined weight space decomposition

\[ S_{\lambda}(V) \cong \oplus_i S_{\lambda}(V)^i .\]

It is "work in progress" to check if Theorem \ref{refined} holds in greater generality and to calculate the character $Char(S_{\lambda}(V)^i)$ for all $\lambda$ and $i$. Such a character formula would imply the classical Weyl character formula.


\begin{example} Flat descent for comodules. \end{example}

In a previous paper I constructed explicit examples of comodules $W_1 \neq W_2$ on $\SLO_{2,\ZZ}$ where the induced comodules $W_1 \otimes \mathbf{Q} \cong W_2 \otimes \mathbf{Q}$ are isomorphic as comodules on 
$\SLO_{2,\mathbf{Q}}$. It is hoped the refined weight space decomposition may be helpful in the study of the problem of classifying such comodules. In general one wants for any integer $d \geq 1$ to classify finite rank torsion free comodules $W$ on $\SL(2,\ZZ)$ such that there is an isomorphism 

\[ W \otimes \mathbf{Q} \cong \Sym^d(V \otimes \mathbf{Q}) \]

of comodules on $\SLO_{2,\mathbf{Q}}$. Here $V$ is the standard right (or left) comodule on $\SLO_{2,\mathbf{Q}}$. The map $\ZZ \rightarrow \mathbf{Q}$ is flat but not faithfully, flat hence we cannot use the classical theory of faithfully flat descent
for comodules to study this problem.

\section{On irreducibility and the refined weight space decomposition}

In this section we prove that the direct summands $\Sym^d(V)^i$ of the symmetric power $\Sym^d(V)$ of the standard comodule $V$ are irreducible comodules on the schematic normalizer $N$ over any field $k$.

Let $k$ be any field and let $G=\Spec(k[G])$ with $k[G]:=k[x_{ij}, 1/D]$ with $D$ the determinant. It follows for any commutative $k$-algebra  $h_G(A)$ is the set of 2 times 2 matrices with coefficients in $A$ and determinant a unit in $A$.
Let $V:=k\{e_1,e_2\}$ be the standard right comodule with the following comodule map:

\[ \Delta_V: V \rightarrow k[G] \otimes_k V \]

defined by 

\[ \Delta_V(e_1):= x_{11} \otimes e_1 + x_{12} \otimes e_2,  \Delta_V(e_2):= x_{21} \otimes e_1 + x_{22} \otimes e_2. \]

Let $k[N]:=R/I$ with variables $t_1:=x_{11}, t_2:=x_{22}, u_1:=x_{12}, u_2:=x_{21}$ as in the previous section.

We get an  induced comodule structure on the $d$'th symmetric power

\[ \Delta: \Sym^d(V) \rightarrow k[N] \otimes \Sym^d(V) \]

defined by

\[\Delta(e_1^ie_2^j):=  t_1^it_2^j \otimes e_1^i e_2^j + u_1^i u_2^j \otimes e_1^j e_2^i.\]

Let $v_i:= e_1^{d-i}e_2^i$ for $i:=0,..,d$. It follows 

\[ \Sym^d(V)^i \cong k\{v_i, v_{d-i} \}.\]

By construction there is for each $i=0,..,d$ an induced comodule map 

\[ \Delta^i. \Sym^d(V)^i \rightarrow k[N] \otimes \Sym^d(V)^i \]

defined as follows:'

\[ \Delta^i(v_i):= t_1^{d-i}t_2^i \otimes v_i + u_1^{d-i}u_2^i \otimes v_{d-i} \]

and 

\[ \Delta^i(v_{d-i}) := t_1^{i}t_2^{d-i} \otimes v_{d-i} + u_1^{i}u_2^{d-i} \otimes v_{i}. \]

\begin{theorem} \label{irreducibleN}  For each $i=0,..,d$ it follows the comodule $(\Sym^d(V)^i, \Delta^i)$ is an irreducible $k[N]$-comodule.
\end{theorem}
\begin{proof} Assume  $v:=a_1v_i + a_2v_{d-i}$ with $a_i \in k$ and let  $L:=kv$. It follows

\[ \Delta^i(v)= \]

\[ (a_1t_1^{d-i}t_2^i, a_2u_1^iu_2^{d-i}) \otimes v_i +  (a_2t_1^{i}t_2^{d-i}, a_1u_1^{d-i}u_2^{i}) \otimes v_{d-i}.\]

It follows $L$ is a submodule iff there is an element $f:=(f_1,f_2) \in k[N]$ with 

\[ \Delta^i(v)= f \otimes v = a_1f \otimes v_i + a_2f \otimes v_{d-i}.\] 

This implies the following equations:

\[ a_1t_1^{d-i}t_2^i =a_1f_1, a_2u_1^iu_2^{d-i}=a_1f_2 \]

and

\[ a_2t_1^it_2^{d-i} = a_2f_1, a_1u_1^{d-i}u_2^i =a_2f_2.\]

If $a_1\neq 0$ it follows $f_1=t_1^{d-i}t_2^i$ and this implies there is an equality

\[ a_2t_1^it_2^{d-i}=a_2f_1= a_2t_1^it_2^{d-i}, \]

hence $a_2=0$. This implies $a_1u_1^{d-i}u_2^i=a_2f_2=0$ hence $a_1=0$ - a contradiction, hence $a_1=0$. This implies
$a_2u_1^iu_2^{d-i}=a_1f_2=0$ hence $a_2=0$. The Theorem follows.
\end{proof}

\begin{example} Irreducible finite dimensional comodules over a field of characteristic $p>0$. \end{example}

If $G:=\GLO(V)$ or $\SLO(V)$ for $V:=k^n$ an $n$-dimensional vector space over $k$, it follows for any finite dimensional irreducible $G$-comodule $V(\lambda)$, there is a Weyl module
$S_{\lambda}(V)$ and an inclusion $V(\lambda) \subseteq S_{\lambda}(V)$ which is strict in general. It is an open problem to calculate the character $Char(V(\lambda))$.
If the terms in the refined weight space decomposition are irreducible comodules over the normalizer $N$, this decomposition may be interesting in the study of such character formulas.
There are canonical diract summand decompositions

\[ V(\lambda) \cong \oplus_i V(\lambda)^i, S_{\lambda}(V) \cong \oplus_j S_{\lambda}(V)^j \]

and if $S_{\lambda}(V)^j$ is an irreducible comodule, it follows the inclusion map 

\[ f_j: V(\lambda)^j \rightarrow S_{\lambda}(V)^j \] 

is either the zero map or an isomorphism.

Note: If the schematic normalizer $N \cong T \times^{\sigma} \Spec(k[W_G])$ is the wreath product of the diagonal torus $T$ with the constand group scheme of the symmetric group $W_G \cong S(n)$ on n letters, we may tro to relate 
comodules on $N$  to comodules on $T$ and $S(n)$. We may try to relate the direct summands in the refined weight space decomposition to comodules on $T$ and $S(n)$.

\section{The symmetric tensors and symmetric powers}

In this section we use the new methods devloped in the previous sections to give explicit examples and applications. Let $(V, \Delta)$ be the standard let comodule on $\GL$ and let $\Sym^2(V), \Sym^2(V)^*, \sym_2(V), \sym_2(V)^*$ be the symmetric powers and symmetric tensors of $V$. We study the weight space decomposition and refined weight space decomposition of the symmetric powers and symmetric tensors and prove they are isomorphic as comodules on $k[N]$ where $T \subseteq \GL$ is the diagonal torus and $T \subseteq N$ is the schematic normalizer of $T$. 



Let $k[G]:=A[x_{ij}, 1/D]$ with $G:=\GL:= \Spec(k[G])$. Let $V:=A\{e_1,e_2\}$ be the free rank $2$ module on the basis $e_1,e_2$ and defined the following comodule map $\Delta$:

\[ \Delta: V \rightarrow V \otimes_A k[G]\]

by

\[ \Delta(e_1):= e_1 \otimes x_{11} + e_2 \otimes x_{21} \]

and

\[ \Delta(e_2):= e_1 \otimes x_{12} + e_2 \otimes x_{22}.\]

\begin{definition} We say the pair $(V, \Delta)$ is the \emph{left standard comodule} on $\GL$.
\end{definition}

Let $C$ be any commutative $A$-algebra and let the tensor power  power $\Sym^2(V)\otimes_A C$ have the basis $U_C:=\{e_1^2 \otimes 1_C, e_1e_2 \otimes1_C, e_2^2\otimes 1_C\}$ where $1_C \in C$ is the multiplicative unit. We get the comodule map

\[ \Delta^2: \Sym^2(V) \otimes_A k[G]  \rightarrow \Sym^2(V)_A \otimes k[G] \]

defined by 

\[  \Delta^2(e_1^2 \otimes 1_C):= e_1^2 \otimes \x^2 + e_1e_2 \otimes (2\x\z) + e_2^2 \otimes \z^2, \]

\[  \Delta^2(e_1e_2 \otimes 1_C):= e_1^2 \otimes \x\y + e_1e_2 \otimes (\x\w + \y \z) + e_2^2 \otimes \z\w, \]

and

\[  \Delta^2(e_2^2 \otimes 1_C):= e_1^2 \otimes \y^2 + e_1e_2 \otimes (2\y\w) + e_2^2 \otimes \w^2 .\]








We may express the map $\Delta^2$ in the basis $U_{k[G]}$ and we get the following formula:

\[ 
M:=[\Delta^2]^{U_{k[G]}}_{U_{k[G]}}=
\begin{pmatrix}    \x^2 &   \x\y &  \y^2 \\
                              2\x\z &  \x\w+ \y\z & 2\y\w \\
                             \z^2 &  \z\w &  \w^2  
\end{pmatrix}.
\]

We may form the tensor power $V^{\otimes 2}:=V \otimes_A V$. Let  $e_{ij}:=e_i \otimes e_j$  and construct  the following comodule map

\[  \Delta^{\otimes 2}: V^{\otimes 2}\otimes k[G]  \rightarrow V^{\otimes 2} \otimes k[G] \]

defined via:

\[ \Delta^{\otimes 2}(e_{11} \otimes 1_{k[G]} ):= e_{11}\otimes  (\x^2) + (e_{12}+ e_{21})\otimes  (\x\z) + e_{22}\otimes  (\z^2) ,\]

\[ \Delta^{\otimes 2}(e_{12} \otimes 1_{k[G]} ):= e_{11}\otimes  (\x\y) + e_{12}\otimes  (\x\w) + e_{21}( \y\z) + e_{22}\otimes  ( \z\w) ,\]

\[ \Delta^{\otimes 2}(e_{21}\otimes 1_{k[G]} ):= e_{11}\otimes   (\x\y) + e_{12}\otimes  (\y\z) + e_{21}\otimes  ( \x\w) + e_{22}\otimes  ( \z\w) ,\]


and

\[ \Delta^{\otimes 2}(e_{22}\otimes 1_{k[G]}  ):= e_{11}\otimes  (\y^2) + (e_{12}+ e_{21}) \otimes  (\y\w) + e_{22}\otimes  (\w^2) .\]

Let $\sym_2(V) \subseteq V\otimes V$ the be invariant under the canonical action of the symmetric group on two elements. Consider $\sym_2(V) \otimes_A C$. Let 

\[  u_{11}:=e_{11} , u_{12}:= e_{12}+e_{21}  , u_{22}:=e_{22}  \]

and let $V_C :=\{ u_{11} \otimes 1_C, u_{12} \otimes 1_C, u_{21} \otimes 1_C$ be the induced basis for $\sym_2(V) \otimes_A C$ as right $C$-module.
It follows $V_{k[G]}$ is a basis for $\sym_2(V) \otimes k[G]$ as right $k[G]$-module.

The comodule map $\Delta^{\otimes 2}$ leaves $\sym_2(V)$ invariant and we get an induced comodule map 

\[ \Delta_2: \sym_2(V)\otimes k[G] \rightarrow \sym_2(V) \otimes k[G] .\]

It follows 

\[ \Delta_2(u_{12}):= u_{11} \otimes(2\x\y) + u_{12} \otimes (\x\w+\y\z) + u_{22} \otimes (2\z\w) .\]

When we express the induced map $\Delta_2$ in the basis $V_{k[G]}$ we get the following matrix:

\[ 
N:=[\Delta_2]^{V_{k[G]}}_{V_{k[G]}}=
\begin{pmatrix}    \x^2 &   2\x\y &  \y^2 \\
                              \x\z &  \x\w+ \y\z & \y\w \\
                             \z^2 &  2\z\w &  \w^2  
\end{pmatrix}.
\]

Let $I:=(x_{12}, x_{21})(x_{11}, x_{22}):=I_1I_2$ be the ideal of the schematic normalizer $N \subseteq G$. It follows there is an isomorphism of commutative rings 
$k[N] \cong A[t_i^m] \oplus A[u_j^n]$ where

\[ \overline{x_{11}}:=t_1, \overline{x_{22}}:=t_2, \overline{x_{12}}:=u_1, \overline{x_{21}}:=u_2,  \]

and where $t_i u_j=0$ in $k[N]$. We get induced comodule maps

\[  \Delta^2: \Sym^2(V) \otimes k[N] \rightarrow \Sym^2(V) \otimes k[N]  \]

and

\[  \Delta_2: \sym_2(V) \otimes k[N] \rightarrow \sym_2(V) \otimes k[N]  \]

If we construct bases $B^{N,2}$ and $C_{N,2}$ for $\Sym^2(V) \otimes k[N]$ and $\sym_2(V) \otimes k[N]$ from the bases $B^2$ and $C_2$  for $\Sym^2(V) \otimes k[G]$ and $\sym_2(V) \otimes k[G]$, 
we get the following matrices:

\[ 
[\Delta^2]^{U_{k[N]}}_{U_{k[N]}}=
\begin{pmatrix}    t_1^2  &   0  &   u_1^2 \\
                              0  &  t_1t_2+u_1u_2  &  0  \\
                             u_2^2 &  0  &  t_2^2   
\end{pmatrix}.
\]

\[ 
[\Delta_2]^{V_{k[N]}}_{V_{k[N]}}=
\begin{pmatrix}    t_1^2  &   0  &   u_1^2 \\
                              0  &  t_1t_2+u_1u_2  &  0  \\
                             u_2^2 &  0  &  t_2^2   
\end{pmatrix}.
\]

\begin{lemma} There is an isomorphism of $k[N]$-comodules $\Sym^2(V) \cong \sym_2(V)$.
\end{lemma}
\begin{proof} This follows from the above calculation since the two matrices  $ [\Delta^2]^{U_{k[N]}}_{U_{k[N]}}  =  [\Delta_2]^{V_{k[N]}}_{V_{k[N]}}$ are "equal", hence they define isomorphic comodules.    \end{proof}

Let us restrict to the diagonal torus $T:=V(I_1)$. It follows $k[T] \cong A[t_i^m]$. We get induce comodules

\[  \Delta^2 : \Sym^2(V) \otimes k[T] \rightarrow \Sym^2(V) \otimes k[T] \]

with corresponding matrix (we let $u_i:=0$)

\[ 
[\Delta^2]^{U_{k[T]}}_{U_{k[T]}}
\begin{pmatrix}    t_1^2  &   0  &    0  \\
                              0  &  t_1t_2  &  0  \\
                              0  &  0  &  t_2^2   
\end{pmatrix}.
\]

We may do a similar construction for $\sym_2(V)$ and we get a map

\[  \Delta_2:  \sym_2(V) \otimes k[T] \rightarrow \sym_2(V) \otimes k[T] \]

with corresponding matrix

\[ 
[\Delta_2]^{V_{k[T]}}_{V_{k[T]}}=
\begin{pmatrix}    t_1^2  &   0  &    0  \\
                              0  &  t_1t_2  &  0  \\
                              0  &  0  &  t_2^2   
\end{pmatrix}.
\]

It follows we can identify the weights of $\Sym^2(V)$ and $\sym_2(V)$: The comodules have the same weights.
Because a map of comodules must preserve the weigth spaces, it follows any map of comodules $\phi: \Sym^2(V)  \rightarrow \sym_2(V)$ must look as follows 

\[ \phi(e_1^2):= au_{11}, \phi(e_1e_2):= bu_{12}, \phi(e_2^2):= cu_{22} \]

with $abc \in A$. A well defined map $\phi$ is an isomorphism iff $abc \in A^*$ is a unit.   If we consider the map $\phi$ as a map 

\[ \phi: \Sym^2(V) \otimes k[G] \rightarrow \sym_2(V) \otimes k[G] \] 

it follows in a choice of bases $U_{k[G]}, V_{k[G]}$, the map $\phi$ will look as follows:

\[ 
U:=[\phi]^{U_{k[G]}}_{V_{k[G]}}=
\begin{pmatrix}     a  &   0  &    0  \\
                              0  &  b  &  0  \\
                              0  &  0  &   c   
\end{pmatrix}.
\]

with inverse map $\phi^{-1}$ given by the matrix 

\[ 
U^{-1}:=[\phi^{-1}]^{U_{k[G]}}_{V_{k[G]}}=
\begin{pmatrix}     1/a  &   0  &    0  \\
                              0  & 1/b  &  0  \\
                              0  &  0  &   1/c   
\end{pmatrix}.
\]

It follows $\Sym^2(V) \cong \sym_2(V)$ iff $U \circ  N \circ U^{-1}=M$ are conjugate matrices. We get the following matrix for $U \circ N \circ U^{-1}$:

\[ 
UNU^{-1}=
\begin{pmatrix}    \x^2 &   2\frac{a}{b}\x\y &  \frac{a}{c}\y^2 \\
                              \frac{b}{a}\x\z &  \x\w+ \y\z & \frac{b}{c}\y\w \\
                             \frac{c}{a}\z^2 &  2\frac{c}{b}\z\w &  \w^2  
\end{pmatrix}.
\]

\begin{theorem}  \label{example} Let $G:=\GL$ with $A$ any commutative ring, and let $N \subseteq G$ be the schematic normalizer of the diagonal torus. There is an isomorphism of $k[N]$-comodules (and $k[T]$-comodules) $\Sym^2(V) \cong \sym_2(V)$.
Let $G:=\GLO_{2,\ZZ}$ with $\ZZ$ the ring of integers. It follows the $k[G]$-comodules $\Sym^2(V)$ and $\sym_2(V)$ are not isomorphic.
\end{theorem} 
\begin{proof}You may check there are no choices of integers $a,b,c \in \ZZ$ with $UNU^{-1}=M$.  Hence an isomorphism $\phi$ of comodules cannot exist.
\end{proof}

\section{Duals of symmetric powers and symmetric tensors}

In this section we study the symmetric powers, symmetric tensors and their duals. We calculate their refined weight space decompositions.

When dualizing a representation of a group scheme $G$  it is better to do this in terms of the functor of points $h_G(-)$ and then to translate this to a dualization functor defined on comodules. At the level of comodules there are several "candidates" for such
a "dualization functor" and one wants such a functor to be related to the representation theory of the group scheme $G$. We want the dual to be functorial and to have the 
property that the double dual $V^{**}$ of a comodule $V$ is canonically isomorphic to $V$. We also want the canonical trace map $tr: V^* \otimes V \rightarrow A$ to be invariant.
The reason one wants such a functor to be defined at the level of comodules is because one wants to calculate the weights of the dual. This is easier done
from the comodule map. We can also calculate the refined weight space decomposition from the comodule map.

Let $k[G]:=A[x_{ij},1/D]$ with $G:=\Spec(k[G]) \cong \GL$. Let $\Delta: W \rightarrow W \otimes_A k[G]$ be a comodule on $G$. Let $\Delta_G$ be the comultiplication on $G$. For any commutative $A$-algebra $B$ 
there is for any $g \in h_G(B)$ a commutative diagram

\[
\diagram      W \rto^{\Delta} \dto    &    W \otimes k[G] \dto^{1\otimes g}   \\
                      W \otimes B \rto^{\Delta(g)}  &   W \otimes B  
\enddiagram
\]

where the map $\Delta(g)$ is defined as follows: Assume $\Delta(w):= \sum_i w_i \otimes y_i$. By definition we let  $\Delta(g)(w \otimes b):= \sum_i w_i \otimes g(y_i)b \in W\otimes B$.

\begin{lemma} \label{multiplicative} The following holds: For any $g,h \in h_G(B)$ it follows $\Delta(g) \circ \Delta(h)= \Delta(g \times h)$ where $g \times h \in h_G(B)$ is the group multiplication. Moreover $\Delta(e)=Id$ if $e \in h_G(B)$ is the identity.
\end{lemma}
\begin{proof} Assume $\Delta(w):= \sum_i w_i \otimes y_i$ and $\Delta(w_i):= \sum_j v(i)_j \otimes x(i)_j$. Assume $\Delta_G(y_i)= \sum_k u(i)_k \otimes w(i)_k$. By definition $g \times h := g \overline{\otimes } h \circ \Delta_G$.

We get the following:

\[  \Delta(g) \circ \Delta(h)(w \otimes b ) := \Delta(g)( \sum_i w_i \otimes h(y_i)b) =\]

\[ \sum_{i,j} v(i)_j \otimes g(x(i)_j)h(y_i) b ) \]

\[  \sum_{i,k} w_i \otimes g(u(i)_k) h(w(i)_k)b  = \]

\[   1 \otimes (g \times h) \Delta(w)b := \Delta(g \times h)(w \otimes b) ,\]

hence $\Delta(g) \circ \Delta(h) = \Delta(g \times h)$. One checks $\Delta(e) = Id$ and the Lemma follows.
\end{proof}

\begin{lemma} \label{canonical} Let $A \rightarrow B$ be a map of commutative rings and let $M,N$ be $A$-modules. There is a canonical map of $B$-modules
\[ \lambda_B: \Hom_A(M,N)\otimes_A B  \rightarrow \Hom_{B}(M \otimes_A B , N \otimes_A B) \]

which is an isomorphism when $M \cong A^n$ and $N:=A$. The map is defined by $\lambda_B(f \otimes b)(m \otimes x):=f(m) \otimes bx$. The map $\lambda_B$ is functorial in $M$.
\end{lemma}
\begin{proof} This is proved in \cite{matsumura}.
\end{proof}

Let $h_G(-)$ be the functor of points of $G$ and define for any comodule $\Delta: W \rightarrow W \otimes k[G]$ the following functor:

\[  W_a(-) : \alg \rightarrow \grp \]

by

\[ W_a(B):= (W \otimes_A B , +)  \]

where we consider $W\otimes_A B$ as an abelian group with addition as group operation. If $\phi: B \rightarrow B'$ is a map of $A$-algebras, there is a canonical map of abelian groups

\[  \phi_a: W_a(B) \rightarrow W_a(B') \]

defined by $\phi_a(w \otimes b):= w \otimes \phi(b)$. Hence $W_a(-)$ is a functor as claimed.

Assume $W \cong A^n$ is a free $A$-module of rank $n$.  There is by Lemma \ref{canonical} for every $B$ a canonical isomorphism

\[ \lambda_B : W^* \otimes_A B \rightarrow (W \otimes_A B)^* .\]

For any element $g \in h_G(B)$ it follows from Lemma \ref{multiplicative} the map $\Delta(g): W \otimes B \rightarrow W \otimes B$ is an isomorphism since $\Delta(g) \circ \Delta(g^{-1})=Id$.  Hence we may define an action

\[ \sigma: h_G(B) \times W^*_A(B) \rightarrow W^*_a(B)  \]

by

\[ \sigma(g, f \otimes b):= \lambda_B^{-1}(\lambda(f \otimes b) \circ \Delta(g^{-1})) .\]

Define the following map: 

\[  \rho_B: h_G(B) \rightarrow \End( W^* \otimes_A B) \]

 by 

\[ \rho_B(g)(f \otimes b):=\sigma(g, f \otimes b).\]

\begin{proposition} \label{dual} For every commutative $A$-algebra $B$ it follows the map $\rho_B$ gives a map of groups 

\[  \rho_B: h_G(B) \rightarrow \operatorname{GL}_B(W^* \otimes_A B).\]

 This construction gives a natural transformation

\[ \Delta^*: h_G(-) \rightarrow h_{\GLO(W^*)}(-)\]

of functors. The representation $(W^*, \Delta^*)$ is the dual of $(W, \Delta)$.
\end{proposition}
\begin{proof} One checks the map $\rho_B$ is a $B$-linear map. If $g,h \in h_G(B)$ it follows

\[  \rho_B(gh)(f \otimes b):= \lambda_B^{-1}(\lambda_B(f \otimes b) \Delta((g \times h)^{-1}) =\]

\[   \lambda_B^{-1}(\lambda_B(f \otimes b)\Delta(h^{-1}) \circ  \Delta(  g^{-1})) =\]

\[ \rho_B(g)( \lambda_B^{-1}(\lambda_B(f \otimes b) \circ \Delta(h^{-1})) := \rho_B(g) \rho_B(h)(f \otimes b) \]

hence $\rho_B(g \times h)= \rho_B(g) \circ \rho_B(h)$. Since $\rho_B(e)$ is the identity map it follows $\rho_B$ is a map of groups. 
The map $\lambda_B$ is functorial and one checks one gets a natural transformation of functors and the Proposition follows.
\end{proof}

\begin{example} The dual of the second symmetric tensor and its weights. \end{example}

We want to calculate the dual representation from Proposition \ref{dual} and the refined weight space decomposition  in some cases. Let $\Delta_1: V \rightarrow V \otimes k[G]$ be the left standard comodule. Recall the following:
$V:=A\{e_1,e_2\}$ is the free rank two module on $e_1,e_2$ and by definition

\[  \Delta_1(e_1):= e_1 \otimes x_{11} + e_2 \otimes x_{21}, \Delta_1(e_2):= e_1 \otimes x_{12} + e_2 \otimes x_{22}. \]

The symmetric tensors is the module of invariants $\sym_2(V):=(V \otimes V)^{S(2)} \subseteq V \otimes V$ under the canonical action of the symmetric group on two elements $S(2)$ acting on $V \otimes V$. Recall the basis $V_C$ for the tensor product
$\sym_2(V) \otimes_A C$.

Define the dual basis of $\sym_2(V)^*$ as the elements 

\[  y_{11}, y_{12}, y_{22}   \]

with $y_{ij}:=u_{ij}^*$. Let $C$ be any commutative $A$-algebra and consider the tensor product  $\sym_2(V)^*\otimes_A C$. It has the following basis

\[  V^*_C:=\{     y_{11} \otimes 1_C , y_{12}\otimes 1_C , y_{22} \otimes 1_C \}.\]

Recall the second symmetric tensor of the comodule map $\Delta_1$ to get a comodule map

\[ \Delta_2: \sym_2(V) \otimes k[G] \rightarrow \sym_2(V) \otimes k[G].  \]

Expressed in the basis $V_C$ we get the following matrix

\[ 
N:=[\Delta_2]^{V_{k[G]}}_{V_{k[G]}}=
\begin{pmatrix}    \x^2 &    2\x\y &    \y^2 \\
                               \x\z &  \x\w+ \y\z &    \y\w \\
                              \z^2 &     2\z\w &  \w^2  
\end{pmatrix}.
\]

The matrix $N$ is invertible with inverse matrix $N^{-1}$ given as follows:

\[ 
N^{-1}=
\begin{pmatrix}    \w^2/D^2  &    -2\y\w/D^2 &    \y^2/D^2 \\
                               -\z\w/D^2 &  (\x\w+ \y\z)/D^2 &    -\x\y/D^2 \\
                              \z^2/D^2 &     -2\x\z/D^2 &  \x^2/D^2  
\end{pmatrix}.
\]

We get an action

\[ \sigma(-,-) : h_G(k[G]) \times \sym_2(V)^*_a(k[G]) \rightarrow \sym_2(V)^*_a(k[G]) \]

defined by 

\[ \sigma(g, f \otimes b):= \lambda_{k[G]}^{-1}(\lambda_{k[G]}(f \otimes b) \circ \Delta(3)(g^{-1})).\]

It follows

\[ \sigma(g, y_{ij} \otimes 1_{k[G]}):= (u_{ij} \otimes 1_{k[G]})^* \circ \Delta(3)(g^{-1}) .\]

Assume $g:=Id_{k[G]} \in h_G(k[G])$. It follows  $g: k[G] \rightarrow k[G ]$ corresponds to the following matrix:

\[ 
g:=
\begin{pmatrix}    \x  &    \y \\
                               \z &  \w  
\end{pmatrix}.
\]

The multiplicative inverse $g^{-1}$ as an element of the group $h_G(k[G])$ corresponds to the matrix

\[ 
g^{-1}:=
\begin{pmatrix}    \w/D  &    -\y/D \\
                               -\z/D &  \x/D  
\end{pmatrix}.
\]

We want to calculate the matrix of the map

\[ \Delta_2(g^{-1}): \sym_2(V) \otimes k[G] \rightarrow \sym_2(V)  \otimes k[G] \]

in the basis $V_{k[G]}$. If $ v\otimes b \in \sym_2(V) \otimes k[G]$ and if $\Delta_2(v):= \sum_i w_i \otimes a_i$, by definition 

\[ \Delta_2(g^{-1})(v \otimes b) := \sum_i w_i\otimes  g^{-1}(a_i)b.\]

It follows

\[ \Delta_2(g^{-1})(u_{11} \otimes b):=\]

\[  u_{11} \otimes (\w^2/D^2)b + u_{12} \otimes (-\z\w/D^2)b + u_{22} \otimes (\z^2/D^2)b , \]
 
\[ \Delta_2(g^{-1})(u_{12} \otimes b):= \]

\[  u_{11} \otimes ((-2\y\w)/D^2)b + u_{12} \otimes ((\x\w+\y\z)/D^2)b + u_{22} \otimes ((-2\x\z) /D^2)b, \]

and

\[ \Delta_2(g^{-1})(u_{22} \otimes b):= \]

\[  u_{11} \otimes (\y^2/D^2)b + u_{12} \otimes (-\x\y/D^2)b + u_{22} \otimes (\x^2/D^2)b. \]

\begin{lemma} There is an equality of matrices 

\[ [\Delta_2(g^{-1})]^{V_{k[G]}}_{V_{k[G]}}=  N^{-1}. \]

\end{lemma}
\begin{proof} The Lemma follows from the above calculation.
\end{proof}

By definition it follows

\[ \lambda_{k[G]}(y_{ij} \otimes 1_{k[G]}):=(u_{ij} \otimes 1_{k[G]})^* .\]

It follows

\[  (u_{11} \otimes 1_{k[G]} )^* \circ \Delta_2(g^{-1})(u_{11} \otimes 1_{k[G]}):= \w^2/D^2 , \]

\[  (u_{11} \otimes 1_{k[G]} )^* \circ \Delta_2(g^{-1})(u_{12} \otimes 1_{k[G]}):= -2\z\w/D^2 , \]

and

\[  (u_{11} \otimes 1_{k[G]} )^* \circ \Delta_2(g^{-1})(u_{22} \otimes 1_{k[G]}):= \y^2/D^2 . \]

It follows

\[ (u_{11} \otimes 1_{k[G]})^* \circ \Delta_2(g^{-1}) = \]

\[  (u_{11} \otimes 1_{k[G]} )(\w^2/D^2)+    (u_{12} \otimes 1_{k[G]} ) (-2\z\w/D^2)+   (u_{22} \otimes 1_{k[G]} )(\y^2/D^2).\]

It follows

\[ \sigma(Id_{k[G]}, y_{11} \otimes 1_{k[G]}):=   \]

\[ y_{11} \otimes (\w^2/D^2)+    y_{12} \otimes (-2\z\w/D^2)+   y_{22} \otimes  (\y^2/D^2),\]

\[ \sigma(Id_{k[G]}, y_{12} \otimes 1_{k[G]}):=   \]

\[ y_{11} \otimes (-\z\w /D^2)+    y_{12} \otimes ((\x\w+\y\z)/D^2)+   y_{22} \otimes  (-\x\y/D^2),\]

and

\[ \sigma(Id_{k[G]}, y_{22} \otimes 1_{k[G]}):=   \]

\[ y_{11} \otimes (\z^2/D^2)+    y_{12} \otimes (-2\x\z/D^2)+   y_{22} \otimes  (\x^2/D^2).\]

We get an induced comodule map

\[ \Delta_2^*: \sym_2(V)^* \otimes k[G] \rightarrow \sym_2(V)^* \otimes k[G]  \]

by defining 

\[ \Delta_2^*(v \otimes b):= \sigma(Id_{k[G]}, v \otimes  b).\]

\begin{lemma}\label{tensormatrix}  There is an equality of matrices

\[ [\Delta_2^*]^{V^*_{k[G]} }_{V^*_{k[G]}} = (N^{-1})^{tr},\]
where $(N^{-1})^{tr}$ is the transpose of $N^{-1}$.
\end{lemma}
\begin{proof} The Lemma follows from the above calculation.
\end{proof}

\begin{example} The "matrix" defining the comodule $\sym_2(V)^*$. \end{example}

Let us write down an explicit formula for the matrix from Lemma \ref{tensormatrix}:

\[ 
(N^{-1})^{tr}=
\begin{pmatrix}    \w^2/D^2  &    -2\z\w/D^2 &    \z^2/D^2 \\
                               -2\y\w/D^2 &  (\x\w+ \y\z)/D^2 &    -2\x\z/D^2 \\
                              \y^2/D^2 &     -\x\y/D^2 &  \x^2/D^2  
\end{pmatrix}.
\]

We get an  induced $k[N]$-comodule map 

\[ \Delta_2^*: \sym_2(V)^* \otimes k[N] \rightarrow \sym_2(V)^* \otimes k[N]. \]

Let $\x:=t_1, \w:=t_2, \y:=u_1, \z:=u_2$ in $k[N]$.  We get an equality of matrices

\[ 
[\Delta_2^*]^{V^*_{k[N]}}_{V^*_{k[N]}}=
\begin{pmatrix}     t_1^{-2}  &     0  &     u_1^{-2} \\
                                0  &   t_1^{-1}t_2^{-1} + u_1^{-1}u_2^{-1}  &    0 \\
                                u_2^{-2}  &    0  &   t_2^{-2}  
\end{pmatrix}.
\]

When we pass to $k[T]$ we get an induced map (also denoted $\Delta_2^*$)

\[ \Delta_2^*: \sym_2(V)^* \otimes k[T] \rightarrow \sym_2(V)^* \otimes k[T] \]

with matrix

\[ 
[\Delta_2^*]^{V^*_{k[T]}}_{V^*_{k[T]}}=
\begin{pmatrix}     t_1^{-2}  &     0  &     0  \\
                                0  &   t_1^{-1}t_2^{-1}   &    0 \\
                                 0   &    0  &   t_2^{-2}  
\end{pmatrix}.
\]

Hence the vectors $\{y_{11}, y_{12}, y_{22} \}$ have weights $\beta_1:=(0,-2),\beta_2:= (-1,-1),\beta_3:= (-2,0)$. Define  $\sym_2(V)^*_{\beta_1, \beta_3}:=A\{y_{11}, y_{22} \}$ and
$\sym_2(V)^*_{\beta_2}:= A\{y_{12} \}$. It follows we get a direct sum decomposition as $k[N]$-comodules

\[  \sym_2(V)^* \cong \sym_2(V)^*_{\beta_1, \beta_3} \oplus \sym_2(V)^*_{\beta_2}. \]

This is the refined weight space decomposition of $\sym_2(V)^*$.

\begin{example} The dual of the symmetric power and its weights. \end{example}

We want to calculate the dual representation from Proposition \ref{dual} and the refined weight space decomposition  in some cases. Let $\Delta_1: V \rightarrow V \otimes k[G]$ be the left standard comodule. Recall the following:
$V:=A\{e_1,e_2\}$ is the free rank two module on $e_1,e_2$ and by definition

\[  \Delta_1(e_1):= e_1 \otimes x_{11} + e_2 \otimes x_{21}, \Delta_1(e_2):= e_1 \otimes x_{12} + e_2 \otimes x_{22}. \]

Consider the symmetric power $\Sym^2(V)$  with basis $B(2):=\{e_1^2 , e_1e_2 , e_2^2\}$ and let $C$ be any commutative $A$-algebra.

For the module $\Sym^2(V) \otimes_A C$ we get the basis $U_R:=\{e_1^2 \otimes 1_R , e_1e_2 \otimes 1_R, e_2^2 \otimes 1_R\}$ where $1_R \in R$ is the multiplicative unit.

We get a  map $\Delta^2$:

\[ \Delta^2 : \Sym^2(V) \otimes k[G] \rightarrow \Sym^2(V) \otimes k[G] .\]

In the basis $U_{k[G]}$ the map $\Delta^2$ is given by the matrix $M$:

\[ 
M:=[\Delta^2]^{U_{k[G]}}_{U_{k[G]}}=
\begin{pmatrix}    \x^2 &    \x\y &    \y^2 \\
                               2\x\z &  \x\w+ \y\z &    \y\w \\
                              \z^2 &     \z\w &  \w^2  
\end{pmatrix}.
\]


We may view the map $M \in \End(\Sym^2(V) \otimes k[G])$ as an invertible endomorphism of $k[G]$-modules.  Define the following map 

\[ 
M^{-1}:=
\begin{pmatrix}    \w^2/D^2 &   - \y\w/D^2 &    \y^2/D^2 \\
                               -2\z\w/D^2 &  (\x\w+ \y\z)/D^2 &    -2\x\y/D^2 \\
                              \z^2/D^2 &     -\x\z/D^2 &  \x^2/D^2  
\end{pmatrix}.
\]

It follows $M^{-1}$ is the multiplicative inverse matrix of $M$. We get an action

\[ \sigma(-,-): h_G(k[G]) \times \Sym^2(V)^*_a(k[G]) \rightarrow \Sym^2(V)^*_a(k[G]) \]

defined by 

\[ \sigma(g, f \otimes b):= \lambda_{k[G ]}^{-1}(\lambda_{k[G]}(f \otimes b) \circ \Delta(g^{-1}).\]

\begin{lemma} If $g:=Id_{k[G]}$   it follows $\Delta^2(g^{-1})=M^{-1}$ is the above matrix. 
\end{lemma}
\begin{proof} The proof is similar to the proof for the symmetric tensor. \end{proof}

Let $\Sym^2(V)^*$ have the basis $B^*:=\{x_1,x_2,x_3\}$ with  $x_1:=(e_1^2)^*, x_2:=(e-1e_2)^*, x_3:=(e_2^2)^*$.
Let $U^*_C:=\{x_1 \otimes 1_C, x_2 \otimes 1_C, x_3 \otimes 1_C\}$. It follows 

\[ \sigma(  Id_{k[G]} , x_i \otimes 1):= \lambda_B^{-1}((v_i \otimes 1)^* \circ M^{-1}).\]  

We get the following calculation:

\[ \sigma(Id_{k[G]}  , x_1 \otimes 1):= x_1 \otimes (\w^2/D^2) + x_2 \otimes (-\y\w/D^2) + x_3\otimes(\y^2/D^2) ,\]

\[ \sigma(Id_{k[G]}  , x_2 \otimes 1):= x_1 \otimes (-2 \z\w/D^2) + x_2 \otimes ((\x\w+\y\z)/D^2) + x_3\otimes(-2 \x\y/D^2) ,\]

and

\[ \sigma(Id_{k[G]}  , x_3 \otimes 1):= x_1 \otimes (\z^2/D^2) + x_2 \otimes (-\x\z/D^2) + x_3\otimes(\x^2/D^2) .\]

We get a map 

\[   \Delta^{2,*}: \Sym^2(V)^* \otimes k[G] \rightarrow \Sym^2(V)^* \otimes k[G] \]

defined as follows: 

Let $U^*_{k[G]}:=\{ x_1 \otimes 1_{k[G]},  x_2 \otimes 1_{k[G]},  x_3 \otimes 1_{k[G]}\}$  be the basis for $\Sym^2(V)^* \otimes k[G]$.
In the basis $U^*_{k[G]}$ we get the matrix

\[ 
(M^{-1})^{tr}:=[\Delta^{2,*}]^{U^*_{k[G]} }_{U^*_{k[G]}}=
\begin{pmatrix}    \w^2/D^2 &   - 2\z\w/D^2 &    \z^2/D^2 \\
                               -\y\w/D^2 &  (\x\w+ \y\z)/D^2 &    -\x\z/D^2 \\
                              \y^2/D^2 &     -2\x\y/D^2 &  \x^2/D^2  
\end{pmatrix}
\]

where $(M^{-1})^{tr}$ is the transpose of $M^{-1}$.

If we pass to the schematic normalizer $k[N]$ with variables $x_{11}:=t_1, x_{22}:=t_2, x_{12}:=u_1, x_{21}:=u_2$ we get the following matrix:

\[ 
[\Delta^{2,*}]^{U^*_{k[N]} }_{U^*_{k[N]}}=   
\begin{pmatrix}      t_1^{-2}  &   0  &    u_1^{-2}  \\
                               0  &   (t_1t_2)^{-1} + (u_1u_2)^{-1}  &   0  \\
                               u_2^{-2}  &     0  &   t_2^{-2}  
\end{pmatrix}.
\]

If we pass to the diagonal torus $T$ (we let $u_i=0$) we get the following:

\[ 
[\Delta^{2,*}]^{U^*_{k[T]} }_{U^*_{k[T]}}=
\begin{pmatrix}      t_1^{-2}  &   0  &    0   \\
                               0  &   t_1^{-1}t_2^{-1}   &   0  \\
                               0   &     0  &   t_2^{-2}   
\end{pmatrix}.
\]

Hence the elements $x_1,x_2,x_3$  have weights $\alpha_1:=(0,-2),\alpha_2:=(-1,-1) ,\alpha_3:=(-2,0)$. Let $\Sym^2(V)_{\alpha_1, \alpha_3 }:=A\{x_1,x_3\}$ and $\Sym^2(V)^*_{\alpha_2}:=A\{x_2\}$.
It follows the direct sum

\[ \Sym^2(V)^* \cong \Sym^2(V)_{\alpha_1, \alpha_3 } \oplus \Sym^2(V)^*_{\alpha_2} \]

is a direct sum of $k[N]$-comodules. It is the refined weight space decomposition of $\Sym^2(V)^*$. We also see directly that $\Sym^2(V)^* \cong \sym_2(V)^*$ are isomorphic 
as $k[N]$-comodules by inpecting the structure matrices - they are "equal". If $A:=\ZZ$ there is no isomorphism $\Sym^2(V)^* \cong \sym_2(V)^*$ of $k[\GLO_{2,\ZZ}]$-comodules since this implies there is (taking double dual)
an isomorphism $\Sym^2(V) \cong \sym_2(V)$ which contradicts Theorem \ref{example}.

\section{The adjoint representation in coordinates}

In this section we give a direct and elementary construction of the adjoint representation of $\GL$  using the coordinate ring and the algebra of distributions. We recover the explicit formula
for the adjoint representation and the refined weight space decompostion constructed in Section 2. 

You should be aware that when using the language of representable functors of Demazure and Gabriel from \cite{DG} 
you are working with the category of all commutative $A$-algebras and this is a non-small category. Hence you are doing "non standard set theory". When defining and proving results using the coordinate ring $k[G]$
you are working with a set $k[G]$ equipped with the structure of an $A$-Hopf algebra. Hence you are doing "classical set theory". You should always understand both viewpoints.


Let $\X,\Y\,\Z,\W$ be independent coordinates over $A$ and let $A[z_{ij}]$
be the polynomial ring in the $z_{ij}$.  Let $G:=\GL$. Define the following actions:

\[ \Delta: A[z_{ij}] \rightarrow A[z_{ij}] \otimes k[G] \]

by 

\[ \Delta(z_{11}) := \X \otimes \x + \Y \otimes \z ,\]

\[ \Delta(z_{12}) := \X \otimes \y + \Y \otimes \w, \]

\[ \Delta(z_{21}) := \Z \otimes \x + \W \otimes \z, \]

and

\[ \Delta(z_{11}) := \Z \otimes \y + \W \otimes \w. \]

Define the map 

\[ \Delta^{-1}: A[z_{ij}] \rightarrow A[z_{ij}] \otimes k[G] \]

by

\[ \Delta^{-1}(z_{11}) := \X \otimes (\w/D) - \Z \otimes (\y/D), \]

\[ \Delta^{-1} (z_{12}) := \Y \otimes (\w/D) - \W \otimes (\y/D), \]  

\[ \Delta^{-1}(z_{21}) := -\X \otimes (\z/D) + \Z \otimes (\x/D), \]

\[ \Delta^{-1}(z_{22}) := -\Y \otimes (\z/D) + \W \otimes (\x/D). \]

\begin{lemma} The pairs $(A[z_{ij}], \Delta), (A[z_{ij}], \Delta^{-1})$ are $k[G]$-comodules.
\end{lemma}
\begin{proof} This follows from an explicit calculation. \end{proof}
 
We may consider the composed map

\[  A[z_{ij}] \rightarrow^{\Delta} A[z_{ij}] \otimes k[G] \rightarrow^{\Delta^{-1} \otimes 1}  A[z_{ij}] \otimes k[G] \otimes k[G] \rightarrow^{m \otimes 1} A[z_{ij}] \otimes k[G] \]

where $m$ is the multiplication map. Let $Ad:= m \otimes 1 \circ \Delta^{-1} \otimes 1 \circ \Delta$.

We get  a map

\[ \Ad: A[z_{ij}] \rightarrow A[z_{ij}] \otimes k[G].\]

\begin{lemma} The pair $(A[z_{ij}], \Ad)$ is a $k[G]$-comodule called the \emph{adjoint representation} of $k[G]$ on $A[z_{ij}]$.
\end{lemma}
\begin{proof} The follows from an explicit calculation. \end{proof}

In coordinates it follows $\Ad$ looks as follows

\[ \Ad(z_{11}) = \X \otimes (\x\w/D) + \Y \otimes (\z\w/D) - \Z \otimes (\x\y/D) - \W \otimes (\y\z/D) , \]

\[ \Ad(z_{12}) = \X \otimes (\y\w/D) + \Y \otimes (\w^2/D) - \Z \otimes (\y^2/D) - \W \otimes (\y\w/D),  \]

\[ \Ad(z_{21}) = - \X \otimes (\x\z/D) - \Y \otimes (\z^2/D) + \Z \otimes (\x^2/D) + \W \otimes (\x\z/D) , \]

and

\[ \Ad(z_{22}) = -  \X \otimes (\y\z/D) - \Y \otimes (\z\w/D) + \Z \otimes (\x\y/D) + \W \otimes (\x\w/D) . \]

Let $I:=(\X-1, \Y, \Z, \W-1) \subseteq A[z_{ij}]$ and let $t_1:=\X-1, e_x:=\Y, e_y:=\Z, t_2:=\W-1$.  It follows

\[ \Ad(t_1) = t_1 \otimes (\x\w/D) + e_x \otimes (\z\w/D) - e_y \otimes (\x\y/D) - t_2 \otimes (\y\z/D) , \]

\[ \Ad(e_x) = t_1 \otimes (\y\w/D) + e_x \otimes (\w^2/D) - e_y \otimes (\y^2/D) - t_2 \otimes (\y\w/D) , \]

\[ \Ad(e_y) = - t_1 \otimes (\x\z/D) - e_x \otimes (\z^2/D) + e_y \otimes (\x^2/D) + t_2 \otimes (\x\z/D) , \]

and

\[ \Ad(t_2) = -  t_1 \otimes (\y\z/D) - e_x \otimes (\z\w/D) + e_y \otimes (\x\y/D) + t_2 \otimes (\x\w/D).  \]

Hence we get an induced map (also denoted $\Ad$):

\[ \Ad: I \rightarrow I \otimes k[G].\]

The map $\Ad$ preserves the square $I^2$ and we get an induced map $\Ad: I/I^2 \rightarrow I/I^2 \otimes k[G]$.

It follows the quotient $I/I^2$ is a free rank $4$ $A$-module on the equivalence classes of the elements $t_1,e_x,e_y,t_2$. I will "abuse the notation" 
and I will not write "overline" to indicate I am calculating with equivalence classes in $I/I^2$. It follows  there is an isomorphism

\[ I/I^2 \cong A\{t_1, e_x, e_y, t_2\}           \]

where the set $\{t_1,e_x,e_y,t_2\}$ is a basis for the $A$-module $I/I^2$. Let $C$ be  any commutative $A$-algebra and consider the tensor product $I/I^2 \otimes_A C$. It has a basis 
$B_C:=\{ t_1 \otimes 1_C, e_x \otimes 1_C, e_y \otimes 1_C, t_2 \otimes 1_C\}$ where $1_C \in C$ is the multiplicative unit. We may view the map $\Ad$ as a map

\[ \Ad: I/I^2 \otimes k[G]  \rightarrow I/I^2 \otimes k[G] \]

and in the basis $B_{k[G]}$ we get the following matrix:

\[ 
M:=[\Ad]^{B_{k[G]}}_{B_{k[G]}}=
\begin{pmatrix}      \x\w/D  &   \y\w/D  &    -\x\z/D &  -\y\z/D   \\
                               \z\w/D   &    \w^2/D   &   -\z^2/D &   -\z\w/D   \\
                               -\x\y/D    &    -\y^2/D   &    \x^2/D &   \x\y/D \\  
                               -\y\z/D &   -\y\w/D &    \x\z/D  &   \x\w/D 
\end{pmatrix}.
\]

The matrix $M$ has the following inverse matrix:

\[ 
M^{-1}=
\begin{pmatrix}      \x\w/D  &   -\x\y/D  &    \y\w/D &  -\y\z/D   \\
                               -\x\z/D   &    \x^2/D   &   -\z^2/D &   \x\z/D   \\
                               \y\w/D    &    -\y^2/D   &    \w^2/D &   -\y\w/D \\  
                               -\y\z/D &   \x\y/D &    -\z\w/D  &   \x\w/D 
\end{pmatrix}.
\]

 Let $T_1:=t_1^*, x:=e_x^*, y:=e_y^*, T_2:=t_2^*$. By Theorem \ref{dist} and Section 2 there is an isomorphism of $A$-modules

\[  \phi: (I/I^2)^* \cong \Dist^+_1(\GL,e) := \Lie(\GL) \cong \gl  \] 

defined by 

\[ \phi(T_1):=e_{11}, \phi(x):=e_{12}, \phi(y):=e_{21}, \phi(T_2):= e_{22}.\]

We get the dual comodule:

\[ \Ad^*: \Lie(\GL) \rightarrow \Lie(\GL) \otimes k[G] \]

defined by the transpose $(M^{-1})^{tr}$ of $M^{-1}$.  It follows $(I/I^2)^* \otimes C$ has the basis:

\[ B^*_C:=\{T_1 \otimes 1_C, x \otimes 1_C, y \otimes 1_C, T_2 \otimes 1_C\}. \]

We get a formula for the matrix of the dual comodule $\Ad^*$:

\[ 
[\Ad^*]^{B^*_{k[G]}}_{B^*_{k[G}}:=(M^{-1})^{tr}=
\begin{pmatrix}      \x\w/D  &   -\x\z/D  &    \y\w/D &  -\y\z/D   \\
                               -\x\y/D   &    \x^2/D   &   -\y^2/D &   \x\y/D   \\
                              \z\w/D    &    -\z^2/D   &    \w^2/D &   -\z\w/D \\  
                               -\y\z/D &   \x\z/D &    -\y\w/D  &   \x\w/D 
\end{pmatrix}.
\]

When we pass to the normalizer $N$ with coordinates $t_i, u_j$ we get an induced comodule

\[ \Ad_N^*: \Lie(\GL) \rightarrow \Lie(\GL) \otimes k[N] \]

defined in the basis $B^*_{k[N]}$ by the following matrix:

\[ 
[\Ad^*_N]^{B^*_{k[N]}}_{B^*_{k[N}}=
\begin{pmatrix}      (1,0)  &    0  &     0  &  (0,1)   \\
                                0    &    t_1t_2^{-1}   &   u_1u_2^{-1} &   0   \\
                               0     &    u_1^{-1}u_2   &     t_1^{-1}t_2  &   0 \\  
                               (0,1) &   0  &    0   &   (1,0)
\end{pmatrix}.
\]

We recover the result we found using the functor of points. We observe there is an equality of  matrices 

\[  [\Ad^*]^{B^*_{k[G]}}_{B^*_{k[G}}=[\Ad]^{B_{k[G]}}_{B_{k[G]}}, \] 

where the last matrix is the one from Section 2 where we defined the adjoint representation using the functor of points. Hence the map $\phi:(I/I^2)^* \rightarrow \gl $ defined by $\phi(T_1):=e_{11}, \phi(x):=e_{12}, \phi(y):=e_{21}, \phi(T_2):=e_{22}$ is an isomorphism of comodules.

When we pass to the diagonal torus we get a formula for the weights of $\Lie(\GL)$:

\[ 
[\Ad^*_T]^{B^*_{k[T]}}_{B^*_{k[T}}=
\begin{pmatrix}      1  &    0  &     0  &  0   \\
                                0    &    t_1t_2^{-1}   &   0  &   0   \\
                               0     &     0   &     t_1^{-1}t_2  &   0 \\  
                               0 &   0  &    0   &   1
\end{pmatrix}.
\]

Hence the vectors $T_1, x,y,T_2$ have weights   $\{(0,0), (-1,1), (1,-1), (0,0)\}$. We also recover the roots of $\Lie(\GL)$ as $R_G:=\{(-1,1), (1,-1)\}$.

We see that it is not so much work to write out the details in terms of the coordinate ring and the comodule map, 
and when we want to calculate the refined weight space decomposition this is easier done in terms of the comodule map. 

You may also compare the above calculation with the calculation of the adjoint representation of $\SL$ on $\sl $ from Section 2.


\end{document}